\documentclass[final,3p]{elsarticle}

\usepackage{amsmath}
\usepackage{amsfonts}
\usepackage{latexsym}

\newtheorem{proposition}{Proposition}[section]

\newtheorem{corollary}{Corollary}[section]
\newtheorem{remark}{Remark}[section]

\journal{ }

\begin{document}

\begin{frontmatter}



\title{Optimization of a plate with holes}


\author[label1,label2]{Dan Tiba}
\ead{dan.tiba@imar.ro}

\author[label3]{Cornel Marius Murea\corref{cor1}}
\ead{cornel.murea@uha.fr}

\address[label1]{Institute of Mathematics (Romanian Academy), Bucharest, Romania}

\address[label2]{Academy of Romanian Scientists, Bucharest, Romania}

\address[label3]{Laboratoire de Math\'ematiques, Informatique et Applications,\\
Universit\'e de Haute Alsace, France}

\cortext[cor1]{Corresponding author}

\begin{abstract}
We consider a simply supported plate with constant thickness, defined 
on an unknown multiply connected domain. We optimize its shape 
according to some given performance functional. Our method is 
of fixed domain type, easy to be implemented, based on a fictitious domain approach and 
the control variational method. The algorithm that we introduce is
of gradient type and performs simultaneous topological and boundary variations.
Numerical experiments are also included and show its efficiency.
\end{abstract}


\end{frontmatter}


\section{Introduction\label{sec:1}}

Shape optimization or optimal design is now a well established 
branch of the calculus of variations. It is a development of the 
optimal control theory with the minimization parameter being just
the domain where the problem is defined. Basic references in this
respect are Pironneau \cite{Pironneau1984}, Sokolowski, Zolesio \cite{Sokolowski1992}, 
Delfour, Zolesio \cite{Delfour2001}, Neittaanm\"aki, Sprekels, Tiba \cite{Tiba2006}, etc.

It is to be noted that the literature on shape optimization problems, 
including unknown or variable domains, is mainly devoted to second order
elliptic equations. Concerning fourth order boundary value problems, 
for instance plate models, there are papers Kawohl, Lang \cite{Kawohl1997}, 
Mu\~noz, Pedregal \cite{Munoz2007}, 
Sprekels, Tiba \cite{Sprekels2003}, 
Arnautu, Langmach, Sprekels, Tiba\cite{Arnautu2000} studying thickness 
optimization problems that may be reduced to optimal control problems by the 
coefficients. In Neittaanm\"aki, Sprekels, Tiba \cite{Tiba2006}, Ch. VI, 
shape optimization problems
for shells and curved rods, with constant thickness, are also studied.
Since their parametric 
representation of the geometric form enters into the coefficients of the model, the shape
optimization problems are again formulated as optimal control problems 
by the coefficients.

It is the aim of this work to extend the study of the optimization and the approximation  
for variable/unknown domain
problems, from the case of second order elliptic operators, to fourth order
operators. The unknowns to be found are the position, the shape, the size,
the number of the holes defining the optimal plate and the given thickness
is assumed constant. The main tools
that we use is the fictitious domain approach 
Neittaanm\"aki, Pennanen, Tiba \cite{Tiba2009}, 
Neittaanm\"aki,  Tiba \cite{Tiba2012},
Halanay, Murea, Tiba \cite{HMT2016}, 
Murea, Tiba \cite{MT2016} 
and the control variational method, 
Barboteu, Sofonea, Tiba \cite{BST2012}, 
Sofonea, Tiba \cite{ST}, 
Neittaanm\"aki, Sprekels, Tiba \cite{Tiba2006} and
their references. 

The plan of the work is as follows. In the next section we discuss the plate model
that we take into account and its approximation via the fictitious domain method,
under weak regularity assumptions on the geometry. This is important from the point
of view of the associated shape optimization problems since it ensures a large class
of admissible domains. Section 3 is devoted to the analysis of such optimal design
problems, including their gradient and a general gradient-type  algorithm,
for their solution. In the last section, numerical examples are investigated
that show the capacity of our approach to generate simultaneous topological
and boundary variations, in the geometric optimization process.
Our results are discussed in $\mathbb{R}^2$ since this is the natural setting
for plates, but extensions to higher dimension are possible.

\section{The model and its approximation}
\setcounter{equation}{0}

Let $\Omega\subset \mathbb{R}^2$ be a bounded, smooth (multiply) connected
open subset representing the shape of a plate of constant thickness
(normalized to one). We consider the fourth order partial differential equation
\begin{eqnarray}
\Delta\Delta y & = & f\hbox{ in }\Omega, \label{2.1}\\
y & = & 0, \quad \Delta y=0\hbox{ on }\partial\Omega, \label{2.2}
\end{eqnarray}
where $f\in L^2(\Omega)$ is the load and $y \in H^4(\Omega)\cap H^1_0(\Omega)$
is the vertical deflection of the plate.
The existence, the regularity and the uniqueness of the strong solution of
(\ref{2.1})-(\ref{2.2}) is well known, under $\mathcal{C}^{1,1}$ conditions
for $\partial\Omega$, \cite{Grisvard1985}.

The difficulty in the numerical solution of (\ref{2.1})-(\ref{2.2}) is that the shape of 
$\Omega$ may be very complicated, if multiply connected, and the standard Finite
Element Method (FEM) may be difficult to implement. Moreover, in the
corresponding shape optimization problems, the geometry may change in each
iteration in a complex way (simultaneous topological and boundary variations)
and this is very costly to be handled by usual discretization methods.

 We consider now another simply
connected smooth bounded domain $D \subset \mathbb{R}^2$ such that $\Omega \subset D$ and
define the following approximation of (\ref{2.1})-(\ref{2.2}), in a sense
to be made precise in the subsequent Proposition \ref{prop:2.1}.
\begin{eqnarray}
-\Delta y_\epsilon + \frac{1}{\epsilon}(1-H_\Omega)y_\epsilon & = & z_\epsilon
\hbox{ in }D, \label{2.3}\\
y_\epsilon & = & 0\hbox{ on }\partial D,\nonumber\\
-\Delta z_\epsilon + \frac{1}{\epsilon}(1-H_\Omega)z_\epsilon & = & f
\hbox{ in }D, \label{2.4}\\
z_\epsilon & = & 0\hbox{ on }\partial D,\nonumber
\end{eqnarray}
where $H_\Omega$ is the characteristic function of $\Omega$ in $D$. 
For the boundary value problems (\ref{2.3})-(\ref{2.4}) we get in the standard
way that the strong solutions satisfy
$y_\epsilon,z_\epsilon\in H^2(D)\cap H^1_0(D)$ if $D$ is in $\mathcal{C}^{1,1}$.
Notice that the systems (\ref{2.3})-(\ref{2.4}) arise from the application
of both the control variational method and fictitious method, as mentioned in
Section \ref{sec:1}.

We relax now the regularity assumptions on the domain $\Omega$ and we suppose that
it is of class $\mathcal{C}$ (the segment property, see \cite{Tiba2006}, \cite{tiba2013} ). In the boundary
value problems  (\ref{2.1})-(\ref{2.2}) and (\ref{2.3})-(\ref{2.4}) we shall
work with weak solutions $y\in H^2(\Omega) \cap H^1_0(\Omega)$ and, respectively, $y_\epsilon, z_\epsilon\in H^1_0(D)$.

\begin{proposition}\label{prop:2.1}
If $\Omega$ is of class $\mathcal{C}$, then $y_\epsilon |_{\Omega}\rightarrow y$
weakly in $H^1_0(\Omega)$ and strongly in $L^2(\Omega)$, where $y \in H^2(\Omega) \cap H_0^1(\Omega)$ satisfies 
(\ref{2.1})-(\ref{2.2}) as a weak solution.
\end{proposition}

\textbf{Proof.} Multiply (\ref{2.4}) by $z_\epsilon$ and integrate by parts:
\begin{equation}\label{2.5}
\int_D |\nabla z_\epsilon|^2 d\mathbf{x} 
+\frac{1}{\epsilon} \int_D (1-H_\Omega)z_\epsilon^2 d\mathbf{x} 
\leq \int_D f\,z_\epsilon\, d\mathbf{x}.
\end{equation}
The Poincar\'e inequality and (\ref{2.5}) gives $\{ z_\epsilon \}$ bounded in
$H^1_0(D)$ and $z_\epsilon \rightarrow \tilde{z}$ strongly in $L^2(D)$
and weakly in $H^1_0(D)$, on a subsequence. Moreover
\begin{equation}\label{2.6}
\int_D (1-H_\Omega)z_\epsilon^2 d\mathbf{x} 
\rightarrow
\int_{D\setminus \Omega} \tilde{z}^2  d\mathbf{x}=0
\end{equation}
due to (\ref{2.5}) since $\{ z_\epsilon \}$ is bounded in $L^p(D)$, $p\geq 1$
in dimension two and we also have $z_\epsilon \rightarrow \tilde{z}$ a.e. in $D$.
One can use Lions' lemma \cite{Lions1969} to infer (\ref{2.6}).
By the Hedberg-Keldys stability property for domains of class $\mathcal{C}$
(see \cite{Tiba2006}) we obtain that $\tilde{z}\in H^1_0(\Omega)$.

The above arguments can be applied to (\ref{2.3}) as well and we have 
$y_\epsilon \rightarrow \tilde{y}$ strongly in $L^2(D)$ and weakly in $H^1_0(D)$
and $\tilde{y}|_\Omega \in H^1_0(\Omega)$. Take any test function 
$\varphi\in \mathcal{C}^\infty_0(\Omega)$ and multiply (\ref{2.3}),
respectively (\ref{2.4}). Since the supports are disjoint, the penalization terms in
(\ref{2.3}), (\ref{2.4}) disappear and we get that $y_\epsilon$ satisfies (\ref{2.1})
in the distribution sense. The boundary condition (\ref{2.2}) are also satisfied due
to the previous remarks. Since the limits $\tilde{y},\ \tilde{z}$ are unique, 
the convergence is in fact valid without taking subsequence.\quad $\Box$

\bigskip

Consider now $H^\epsilon : D \rightarrow \mathbb{R}$ to be a $C^1$ regularization of the 
characteristic function $H_\Omega$ and $H^\epsilon\rightarrow H_\Omega$ strongly
in $L^p(\Omega)$, $p\geq 1$. Examples of this type will be indicated
in the next section.

\begin{corollary}\label{cor:2.2}
If in (\ref{2.3}), (\ref{2.4}) we replace $H_\Omega$ by $H^\epsilon$, the other notations
being preserved, then the conclusion of Proposition \ref{prop:2.1} remains valid.
\end{corollary}

\section{Shape optimization problems and their gradient}
\setcounter{equation}{0}

We associate to (\ref{2.1}), (\ref{2.2}) the following minimization problem
\begin{equation}\label{3.1}
\min_{\Omega\in \mathcal{O}}
\int_{\Lambda} J\left(\mathbf{x},y(\mathbf{x})\right) d\mathbf{x},
\end{equation}
where $\mathcal{O}$ is the class of admissible domains to be defined below,
$y\in H^1_0(\Omega)$ is the weak solution of (\ref{2.1}), (\ref{2.2}), 
$\Lambda$ may be $\Omega$ or $\partial\Omega$ or some part of $\Omega$ or $\partial\Omega$
and $J$ is the performance index of Carath\'eodory type (measurable in $\mathbf{x}$
and continuous in $y$). More hypotheses or constraints will be imposed as necessity
appears. The problem (\ref{3.1}), (\ref{2.1}), (\ref{2.2}) has a similar form with
optimal control problems, however the optimization parameter here is the 
geometry, the domain $\Omega$ itself.

The family $\mathcal{O}$ should be ``large'' in order to perform the
optimization in (\ref{3.1}) on a consistent admissible class. 
We avoid regularity hypotheses on the geometry (that are frequently used in shape
optimization, see \cite{Chenais1975}, \cite{Pironneau1984}, \cite{Sokolowski1992})
and we have just assumed that any $\Omega\in \mathcal{O}$ is 
an open set of class $\mathcal{C}$, contained in some given bounded domain
$D\subset\mathbb{R}^2$:
\begin{equation}\label{3.2}
\Omega \subset D,\quad \forall \Omega\in \mathcal{O}.
\end{equation}
On may add the constraint
\begin{equation}\label{3.3}
E \subset \Omega,\quad \forall \Omega\in \mathcal{O}
\end{equation}
where $E \subset\subset D$ is some given not empty subset of $\mathbb{R}^2$.

Let $X(D)$ denote a subset of $\mathcal{C}(\overline{D})$. For instance, $X(D)$ may
be a finite element space defined in $D$.
Following \cite{Tiba2009}, \cite{Tiba2012}, with any $g\in X(D)$, that we
call a parametrization of the geometry, we associate the open set
\begin{equation}\label{3.4}
\Omega_g = int \left\{\mathbf{x} \in D;\ g(\mathbf{x})\geq 0\right\}.
\end{equation}
In the absence of regularity assumptions and due to the possible presence of
critical points  of $g$, it is possible that
$g$ has level set $\{ \mathbf{x} \in D;\ g(\mathbf{x}) = k \}$ of positive measure.
This is the reason for the form of the definition (\ref{3.4}). This is, in principle,
different from the set of points where $g(x) > 0$. Notice that $\Omega_g$ is
a Carath\'eodory
open set, i.e. cracks or cuts are not allowed. However, high oscillations
of the boundary
are possible (and the segment property may not be always valid and has to
be imposed separately).
In general, $\Omega_g$ may have   many connected components, that may be
multiply connected.
If constraint (\ref{3.3}) is imposed, then $X(D)$ should include the condition:
\begin{equation}\label{3.5}
g(\mathbf{x})\geq 0\hbox{ in }E.
\end{equation}

If $H:\mathbb{R}\rightarrow \mathbb{R}$ denotes the maximal monotone extension of the
Heaviside function (see \cite{Tiba2009}, \cite{1992}) then $H(g)$ is the
characteristic function of $\overline{\Omega}_g$.
The regularization $H^\epsilon=H^\epsilon(g)$, from Corollary \ref{cor:2.2}, can be
simply obtained
by a regularization of the Heaviside function. In \cite{Tiba2012}, the following
formula is used
\begin{equation}\label{3.6}
H^\epsilon(r)=
\left\{
\begin{array}{ll}
1-\frac{1}{2}e^{-\frac{r}{\epsilon}}, & r \geq 0,\\
\frac{1}{2}e^{\frac{r}{\epsilon}}, & r < 0
\end{array}
\right.
\end{equation}
but other choices are possible.

Taking into account the approximation results from the previous section, we
approximate the
minimization problem (\ref{3.1}), (\ref{2.1}), (\ref{2.2}) by (\ref{3.1}), (\ref{2.3}), (\ref{2.4})
where $H_\Omega$ is replaced by $H^\epsilon(g)$. The cost functional (\ref{3.1}),
depending on the form of $\Lambda$, may be approximated in the form
\begin{equation}\label{3.7}
\int_E J\left(\mathbf{x},y_\epsilon(\mathbf{x})\right) d\mathbf{x},
\quad\hbox{if }\Lambda=E,
\end{equation}
\begin{equation}\label{3.8}
\int_D H^\epsilon(g) J\left(\mathbf{x},y_\epsilon(\mathbf{x})\right) d\mathbf{x},
\quad\hbox{if }\Lambda=\Omega .
\end{equation}

The case $\Lambda=\partial\Omega$ imposes more regularity assumptions on the geometry
in order to ensure the application of trace theorems and it has been recently discussed in 
Tiba \cite{2017} for second order operators. We limit our investigations here to (\ref{3.7}), (\ref{3.8}).
The approximation of the state equation and of the cost functionals ensures that all the 
computations are to be performed in the fixed domains $E$ or $D$.
The geometry $\Omega$ is hidden under this approach in the mapping $g\in X(D)$. Consequently,
the approximating shape optimization problems are in fact optimal control problems with
the control $g$ acting in the coefficients of the lowest order term in the differential operator.
Notice as well the smooth dependence of $y_\epsilon$ on  $g$ when $H^\epsilon$ is used instead of $H$.
This is analyzed in the next result and is fundamental for the application of the gradient
methods in the solution of the optimization problem (\ref{3.1}), (\ref{2.1}), (\ref{2.2}).

\begin{proposition}\label{prop:3.1}
The mappings $g\rightarrow y_\epsilon=y_\epsilon(g)$, $g\rightarrow z_\epsilon=z_\epsilon(g)$
defined by (\ref{2.3}), (\ref{2.4}) with $H_\Omega$ replaced by $H^\epsilon(g)$ are
G\^ateaux differentiable between $\mathcal{C}(D)$ and $H^1_0(\Omega)$ and
$w=\nabla y_\epsilon(g)v$, $u=\nabla z_\epsilon(g)v$ for any $v$ in $\mathcal{C}(D)$
satisfy the following system in variations:
\begin{eqnarray*}
-\Delta u + \frac{1}{\epsilon}(1-H^\epsilon(g)) u & = & 
\frac{1}{\epsilon} (H^\epsilon)^\prime(g) z_\epsilon v,\\
-\Delta w + \frac{1}{\epsilon}(1-H^\epsilon(g)) w & = & 
u+ \frac{1}{\epsilon} (H^\epsilon)^\prime (g) y_\epsilon v, \\
\end{eqnarray*}
with $u,w\in H^1_0(\Omega)$.
\end{proposition}

\textbf{Proof.} We denote by $y_\epsilon^\lambda=y_\epsilon(g+\lambda v)$,
$z_\epsilon^\lambda=z_\epsilon(g+\lambda v)$, $\lambda\in \mathbb{R}$.
Substrating the corresponding regularized equations and dividing by $\lambda\neq 0$, we get
\begin{eqnarray}
-\Delta \frac{z_\epsilon^\lambda-z_\epsilon}{\lambda} 
+ \frac{1}{\epsilon}(1-H^\epsilon(g+\lambda v)) \frac{z_\epsilon^\lambda-z_\epsilon}{\lambda} & = & 
\frac{1}{\epsilon} \frac{H^\epsilon(g+\lambda v)-H^\epsilon(g)}{\lambda} z_\epsilon ,\label{3.9}\\
-\Delta \frac{y_\epsilon^\lambda-y_\epsilon}{\lambda} 
+ \frac{1}{\epsilon}(1-H^\epsilon(g+\lambda v)) \frac{y_\epsilon^\lambda-y_\epsilon}{\lambda} & = & 
 \frac{z_\epsilon^\lambda-z_\epsilon}{\lambda}\nonumber\\
&& +\frac{1}{\epsilon} \frac{H^\epsilon(g+\lambda v)-H^\epsilon(g)}{\lambda} y_\epsilon , \label{3.10}
\end{eqnarray}
will null boundary conditions on $\partial D$ for 
$y_\epsilon, z_\epsilon, y_\epsilon^\lambda,z_\epsilon^\lambda$. Here $\epsilon>0$ is fixed and 
$\lambda\in \mathbb{R}$ is the varying parameter ($\lambda\rightarrow 0$).

We multiply (\ref{3.9}) by $\frac{z_\epsilon^\lambda-z_\epsilon}{\lambda}$ and, after some computations, 
we  get
\begin{eqnarray}
 \int_D \left| \nabla \frac{z_\epsilon^\lambda-z_\epsilon}{\lambda} \right|^2 d\mathbf{x}
+ \frac{1}{\epsilon}\int_D (1-H^\epsilon(g+\lambda v)) 
\left| \frac{z_\epsilon^\lambda-z_\epsilon}{\lambda} \right|^2 d\mathbf{x}\nonumber\\
=  
\frac{1}{\epsilon} \int_D 
\frac{H^\epsilon(g+\lambda v)-H^\epsilon(g)}{\lambda} 
z_\epsilon \frac{z_\epsilon^\lambda-z_\epsilon}{\lambda}d\mathbf{x} . \label{3.11}
\end{eqnarray}

Since $H^\epsilon$ is of class $\mathcal{C}^1$, we have 
$\frac{H^\epsilon(g+\lambda v)-H^\epsilon(g)}{\lambda} \rightarrow (H^\epsilon)^\prime (g) v$ a.e. in $D$
and it is bounded in $L^\infty(D)$ with respect to $\lambda\in \mathbb{R}$. We get from (\ref{3.11})
that $\left\{ \frac{z_\epsilon^\lambda-z_\epsilon}{\lambda}\right\}$ is bounded in $H^1_0(D)$.
On a subsequence, we have $\frac{z_\epsilon^\lambda-z_\epsilon}{\lambda} \rightarrow u \in H^1_0(D)$,
weakly in $H^1_0(D)$ and strongly in $L^2(D)$.

A similar argument, using the boundedness of $\left\{ \frac{z_\epsilon^\lambda-z_\epsilon}{\lambda}\right\}$
applied to (\ref{3.10}), gives that 
$\left\{ \frac{y_\epsilon^\lambda-y_\epsilon}{\lambda}\right\}$
is bounded in $H^1_0(D)$ and converges weakly in $H^1_0(D)$ and strongly in $L^2(D)$,
to some limit $w \in H^1_0(D)$, on a subsequence.
Passing to the limit in (\ref{3.9}), (\ref{3.10}) on a common subsequence, we get the equations from the
proposition, satisfied by $u,w \in H^1_0(D)$.

We notice that the equations for $u$, $w$ have a unique solution and this shows that the above
convergences are valid without taking subsequences. We conclude the G\^ateaux differentiability of the
maps $y_\epsilon(g)$,  $z_\epsilon(g)$ and the proof is finished.\quad $\Box$

\medskip

We introduce now the so called adjoint system. To do this, we shall consider two cases of the cost
functionals:
\begin{equation}\label{3.12}
\frac{1}{2}\int_E (y_\epsilon-y_d)^2 d\mathbf{x},
\end{equation}
which is a special case of (\ref{3.7}) with some given $y_d \in L^2(D)$. The second functional is 
(\ref{3.8}).

For the performance index (\ref{3.12}), we introduce the following adjoint system
\begin{eqnarray}
-\Delta p + \frac{1}{\epsilon}(1-H^\epsilon(g)) p & = & 
\chi_E (y_\epsilon-y_d)\hbox{ in }D,
\label{3.13}
\\
-\Delta q + \frac{1}{\epsilon}(1-H^\epsilon(g)) q & = & 
p \hbox{ in }D, \label{3.14}
\\
p=0,\quad q&=&0 \hbox{ on }\partial D, \label{3.15}
\end{eqnarray}
where $\chi_E$ is the characteristic function of $E$ in $D$.

\begin{proposition}\label{prop:3.2}
The directional derivative of the cost functional (\ref{3.12}) is given by
$$
\frac{1}{\epsilon}\int_D (H^\epsilon)^\prime (g) v (y_\epsilon p +z_\epsilon q)d\mathbf{x},
$$
for $p,q$ satisfying (\ref{3.13})--(\ref{3.15}) and for any $v\in \mathcal{C}(D)$.
\end{proposition}

\textbf{Proof.} We have (in the notations of Proposition \ref{prop:3.1}):
\begin{eqnarray*}
L &= & \lim_{\lambda\rightarrow 0}\frac{1}{2\lambda}\left[
\int_E (y_\epsilon^\lambda-y_d)^2 d\mathbf{x} - \int_E (y_\epsilon-y_d)^2 d\mathbf{x}
\right]
= \lim_{\lambda\rightarrow 0} \int_E \frac{y_\epsilon^\lambda-y_\epsilon}{\lambda}
\, \frac{y_\epsilon^\lambda + y_\epsilon -2y_d}{2} d\mathbf{x}\\
 &= & \int_E w (y_\epsilon-y_d) d\mathbf{x}
=\int_D  w \left(-\Delta p + \frac{1}{\epsilon}(1-H^\epsilon(g)) p\right) d\mathbf{x}\\
 &= & \int_D p \left(-\Delta w + \frac{1}{\epsilon}(1-H^\epsilon(g)) w\right) d\mathbf{x},
\end{eqnarray*}
by (\ref{3.13}) and the partial integration.

Using Proposition \ref{prop:3.1}, we get:
\begin{eqnarray*}
L &= & \int_D p \left(u+ \frac{1}{\epsilon} (H^\epsilon)^\prime (g) y_\epsilon v\right) d\mathbf{x}\\
&= & \frac{1}{\epsilon}\int_D (H^\epsilon)^\prime (g) y_\epsilon p\, v \, d\mathbf{x}
+\int_D \left(-\Delta q + \frac{1}{\epsilon}(1-H^\epsilon(g)) q\right) u \, d\mathbf{x}\\
&= & \frac{1}{\epsilon}\int_D (H^\epsilon)^\prime (g) y_\epsilon p\, v \, d\mathbf{x}
+\int_D \left(-\Delta u + \frac{1}{\epsilon}(1-H^\epsilon(g)) u\right) q \, d\mathbf{x}\\
&= & \frac{1}{\epsilon}\int_D (H^\epsilon)^\prime (g) y_\epsilon p\, v \, d\mathbf{x}
+\frac{1}{\epsilon}\int_D (H^\epsilon)^\prime (g) z_\epsilon q\, v \, d\mathbf{x}
\end{eqnarray*}
by (\ref{3.14}) and again by Proposition \ref{prop:3.1}. This ends the proof.\quad $\Box$

\medskip

If the cost functional (\ref{3.8}) is taken into account, the equation in variation
is given by Proposition \ref{prop:3.1} as well, but in the adjoint system (\ref{3.13})--(\ref{3.15}),
the equation (\ref{3.13}) has to be replaced by
\begin{equation}\label{3.16}
-\Delta p + \frac{1}{\epsilon}(1-H^\epsilon(g)) p  =  
H^\epsilon(g) J_y^\prime(\mathbf{x},y_\epsilon)v\hbox{ in }D,
\end{equation}
under the differentiability assumption for $J(\mathbf{x},\cdot)$ and the $L^2(D)$
integrability for $J_y^\prime(\mathbf{x},y_\epsilon)$. In a similar way, we get

\begin{corollary}\label{cor:3.1}
The directional derivative of the cost functional (\ref{3.8}) has the form:
$$
\int_D (H^\epsilon)^\prime (g)\left[ 
J(\mathbf{x},y_\epsilon(\mathbf{x}))
+\frac{1}{\epsilon}\left(y_\epsilon(\mathbf{x}) p(\mathbf{x}) +z_\epsilon(\mathbf{x}) q(\mathbf{x}) \right)
\right]
v(\mathbf{x})\, d\mathbf{x}.
$$
\end{corollary}

The first term in the above formula appears since in (\ref{3.8}) the derivative of $H^\epsilon(g)$,
for perturbation $g+\lambda v$, also appears.

\begin{remark}\label{rem:3.1}
By Proposition \ref{3.2}, the gradient of the performance index (\ref{3.12})
is 
$\frac{1}{\epsilon} (H^\epsilon)^\prime (g) (y_\epsilon p +z_\epsilon q)$ 
and the steepest descent direction is  with minus sign.
Another descent direction is $-(y_\epsilon p +z_\epsilon q)$ since the coefficient is positive 
due to the monotocity of $H^\epsilon(\cdot)$. It also has the advantage of simplicity.
If polynomial regularizations of $H_\Omega$, like
$$
\widetilde{H}^\epsilon(r)=
\left\{
\begin{array}{ll}
1, & r \geq 0,\\
\frac{\epsilon(r+\epsilon)^2-2r(r+\epsilon)^2}{\epsilon^3}, & -\epsilon < r <0,\\
0, & r \leq -\epsilon
\end{array}
\right.
$$
are used instead of (\ref{3.6}), then the support of the gradient or 
of the steepest descent direction is in the set $\{ -\epsilon < g(\mathbf{x}) <0 \}$,
that is in a neighborhood of $\partial\Omega_g$ (when the roots of $g(\cdot)$ are
noncritical).
Similar considerations may be made in connection to the functional (\ref{3.8})
and Corollary \ref{cor:3.1}. Both variants of descent directions may generate boundary
and/or topological variations of the domain $\Omega_g$. A more general situation is considered in the Proposition \ref{prop:4.1}, in the next section.
\end{remark}

As we have already mentioned, the shape optimization problem
(\ref{3.1}), (\ref{2.1}), (\ref{2.2}) may be approximated by (\ref{3.1}), (\ref{2.3}), (\ref{2.4}).
Using admissible domains defined in (\ref{3.4}) and regularizations like (\ref{3.6})
with approximation of the characteristic functions 
$H_{\Omega_g}$ by $H^\epsilon(g)$, we have to solve an optimal control problem
with control $g\in X(D)$ acting in the lower order terms of the system.
In particular, we also infer the necessary optimality conditions for the approximating control 
problem (\ref{3.1}), (\ref{2.1}), (\ref{2.2}) with $H^\epsilon(g)$ instead of $H_{\Omega_g}$.

\begin{corollary}\label{cor:3.2}
Let $g_\epsilon^*\in X(D)$ denote an optimal solution. The optimality
conditions for $g_\epsilon^*$ are given by the system (\ref{3.1}), (\ref{2.3}), (\ref{2.4}),
the adjoint system (\ref{3.13})--(\ref{3.15}) (or (\ref{3.14})--(\ref{3.16}) according
to the form (\ref{3.12}), respectively (\ref{3.8}) of the cost) and the maximum 
principle:
$$
\int_D (H^\epsilon)^\prime (g_\epsilon^*)(y_\epsilon^* p_\epsilon^* 
+z_\epsilon^* q_\epsilon^*)v\, d\mathbf{x} \leq 0,
\quad \forall v,
$$
respectively
$$
\int_D (H^\epsilon)^\prime (g_\epsilon^*)\left[
J(\mathbf{x},y_\epsilon^*(\mathbf{x}))
+\frac{1}{\epsilon}
(y_\epsilon^*(\mathbf{x}) p_\epsilon^*(\mathbf{x}) +z_\epsilon^*(\mathbf{x}) q_\epsilon^*(\mathbf{x}))
\right] v(\mathbf{x})\, d\mathbf{x} \leq 0,
\quad \forall v,
$$
where $y_\epsilon^*,z_\epsilon^*\in H^1_0(D)$ denote the approximating optimal states,
$p_\epsilon^*, q_\epsilon^*$ denote the corresponding adjoint states and 
$v\in \mathcal{C}(\overline{D})$ is any admissible variation such that 
$g_\epsilon^*+\lambda v \in X(D)$ for $\lambda>0$, small.

For instance, if  $X(D)$ is given by (\ref{3.5}), the admissible $v$ have to satisfy
(\ref{3.5}) as well.
\end{corollary}

By Proposition \ref{prop:3.2} and Corollary \ref{cor:3.1}, gradient methods may be applied
with various descent directions. We formulate the following general gradient with projection
algorithm:

\newpage

\textbf{Algorithm 3.1}

\vspace{2mm}

\textbf{Step 1} Start with $n=0$, $\epsilon >0$ given ``small'' and select some initial $g_n$.

\textbf{Step 2} Compute $y_\epsilon^n,\ z_\epsilon^n$ the solution of (\ref{2.3}), (\ref{2.4})
with $H_{\Omega_g}$ replaced by $H^\epsilon(g)$.

\textbf{Step 3} Compute $p^n,\ q^n$ the solution of (\ref{3.13})--(\ref{3.15}) or
(\ref{3.14})--(\ref{3.16}).

\textbf{Step 4} Compute the gradient of the considered cost functional according to 
Proposition \ref{prop:3.2}, respectively Corollary \ref{cor:3.1}.

\textbf{Step 5} Denote by $w_n$ the chosen descent direction, according to Remark 3.1
and define $\widetilde{g}_n=g_n+\lambda_n w_n$, where $\lambda_n >0$ is obtained via some line search.

\textbf{Step 6} Compute $g_{n+1}=Proj_{X(D)} (\widetilde{g}_n)$, if the constraint (\ref{3.5}) is imposed.

\textbf{Step 7} If $| g_n-g_{n+1}|$  and/or $| \nabla j (g_n)|$ are below some prescribed 
tolerance parameter, then \textbf{Stop}.
If not, update $n:=n+1$ and go to \textbf{Step 2}.

\vspace{1mm}
Notice that, according to \cite{ptiba}, in case constraints are imposed on $g$ (for instance, as in \textbf{Step 6}), the set $X(D)$ should consist of piecewise continuous functions, due to the projection operation. The above arguments can be extended to this case in a rather straightforward way.

\noindent
In all the examples discussed in the next section, we underline the combination of both topological and boundary variations that is a property of Algorithm 3.1.
\bigskip

\section{Numerical examples}
\setcounter{equation}{0}

We have employed the software FreeFem++, \cite{freefem++}.
 
\medskip

\textbf{Example 1.} 

This is inspired by the example 2 from \cite{Tiba2009}, but the second order elliptic
equation is replaced by (\ref{2.1})--(\ref{2.2}).
We have $D=]-1,1[\times ]-1,1[$, the 
load $f=3$, the cost function $j(g)=\frac{1}{2}\int_\Omega (y_\epsilon-y_d)^2 d\mathbf{x}$, where
$y_d(x_1,x_2)=-(x_1-0.5)^2 -(x_2-0.5)^2+\frac{1}{16}$.
The initial geometric parametrization function is 
$$
g_0(x_1,x_2)=\min \left(x_1^2+x_2^2-\frac{1}{16};\ (x_1-0.5)^2 + x_2^2-\frac{1}{64};
\ 1-x_1^2-x_2^2 \right),
$$ 
which corresponds to a domain with two holes (see Fig.1).
We use for $D$ a mesh of 53360 triangles and 26981 vertices and 
for the approximation of $g$, $y$, $z$ we use piecewise linear 
finite elements, globally continuous (no constraints on $g$). 
The penalization parameter is $\epsilon=10^{-5}$. 

From Corollary \ref{cor:3.1} and Remark \ref{rem:3.1}, we get that
$$
-\left[ 
J(\mathbf{x},y_\epsilon(\mathbf{x}))
+\frac{1}{\epsilon}\left(y_\epsilon(\mathbf{x}) p(\mathbf{x}) +z_\epsilon(\mathbf{x}) q(\mathbf{x}) \right)
\right]
$$
is a descent direction. 
The cost functional is of type (\ref{3.8}) with 
$J\left(\mathbf{x},y_\epsilon(\mathbf{x})\right)=\frac{1}{2}(y_\epsilon-y_d)^2$.
First, accroding to Algorithm 3.1, we use the descent direction
\begin{equation}\label{4.2}
w_n=-\left[ 
\frac{1}{2}(y_\epsilon^n-y_d)^2
+\frac{1}{\epsilon}(y_\epsilon^n p^n + z_\epsilon^n q^n)
\right].
\end{equation}

The sequence $\left( j(g_{n}) \right)_{n\in\mathbb{N}}$ is decreasing.
For the stopping test, we can use: \texttt{if} $ | j(g_{n}) | < tol$ 
\texttt{then} STOP, where $tol=10^{-10}$. To simplify the notation, we write $\Omega_n$ in place of $\Omega_{g_n}$.

\begin{figure}[ht]
\begin{center}
\includegraphics[width=5cm]{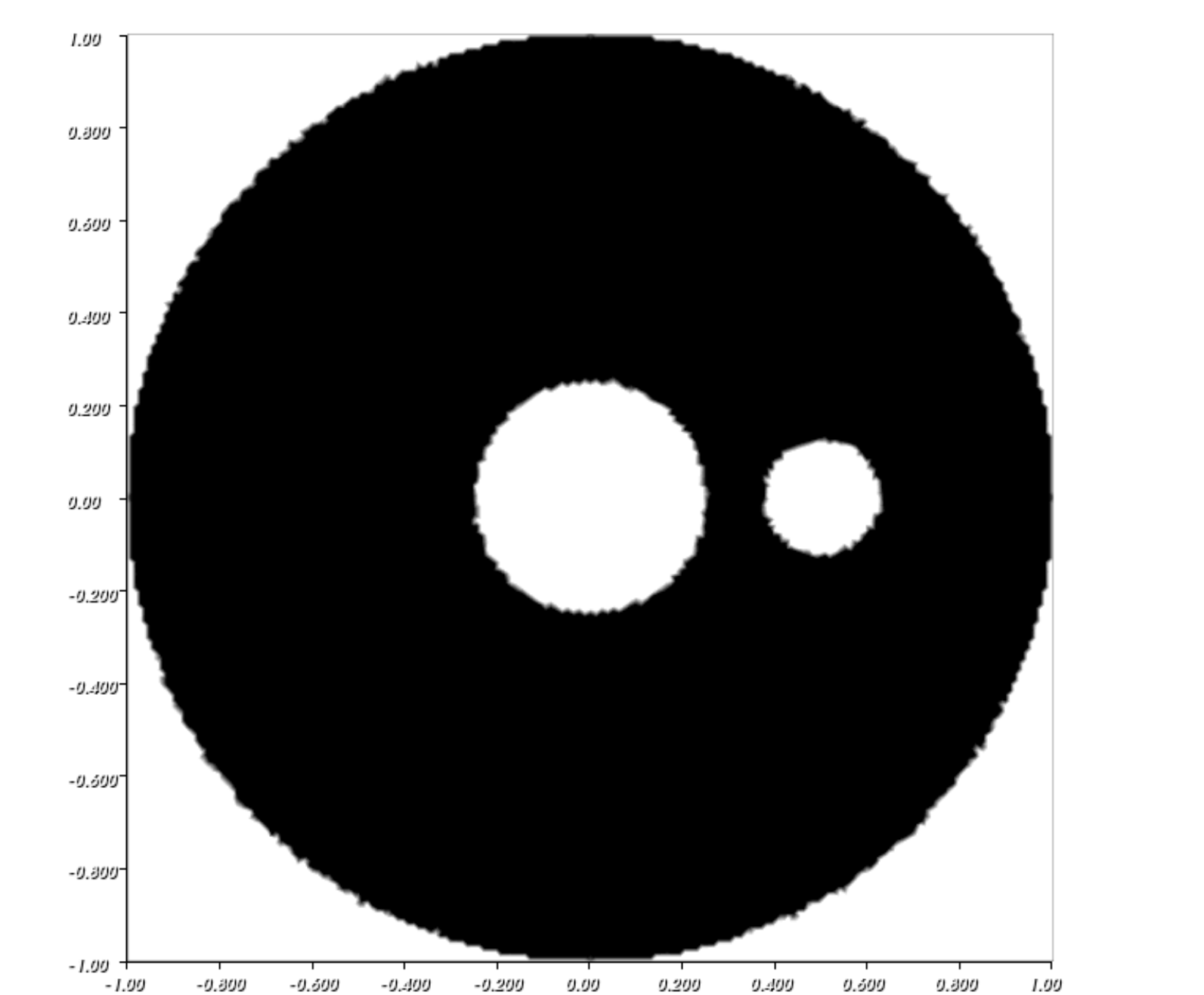}
\quad
\includegraphics[width=5cm]{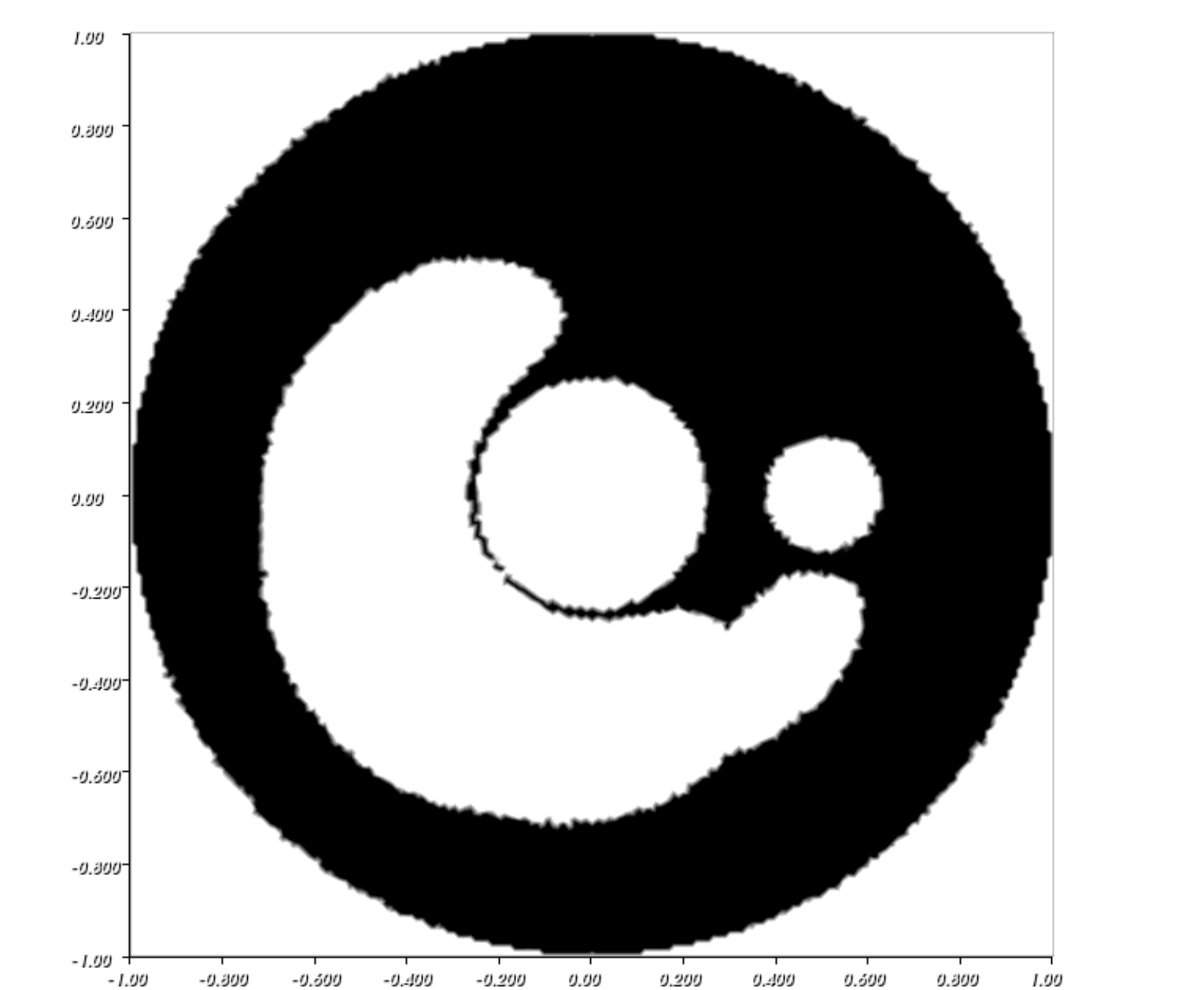}
\quad
\includegraphics[width=5cm]{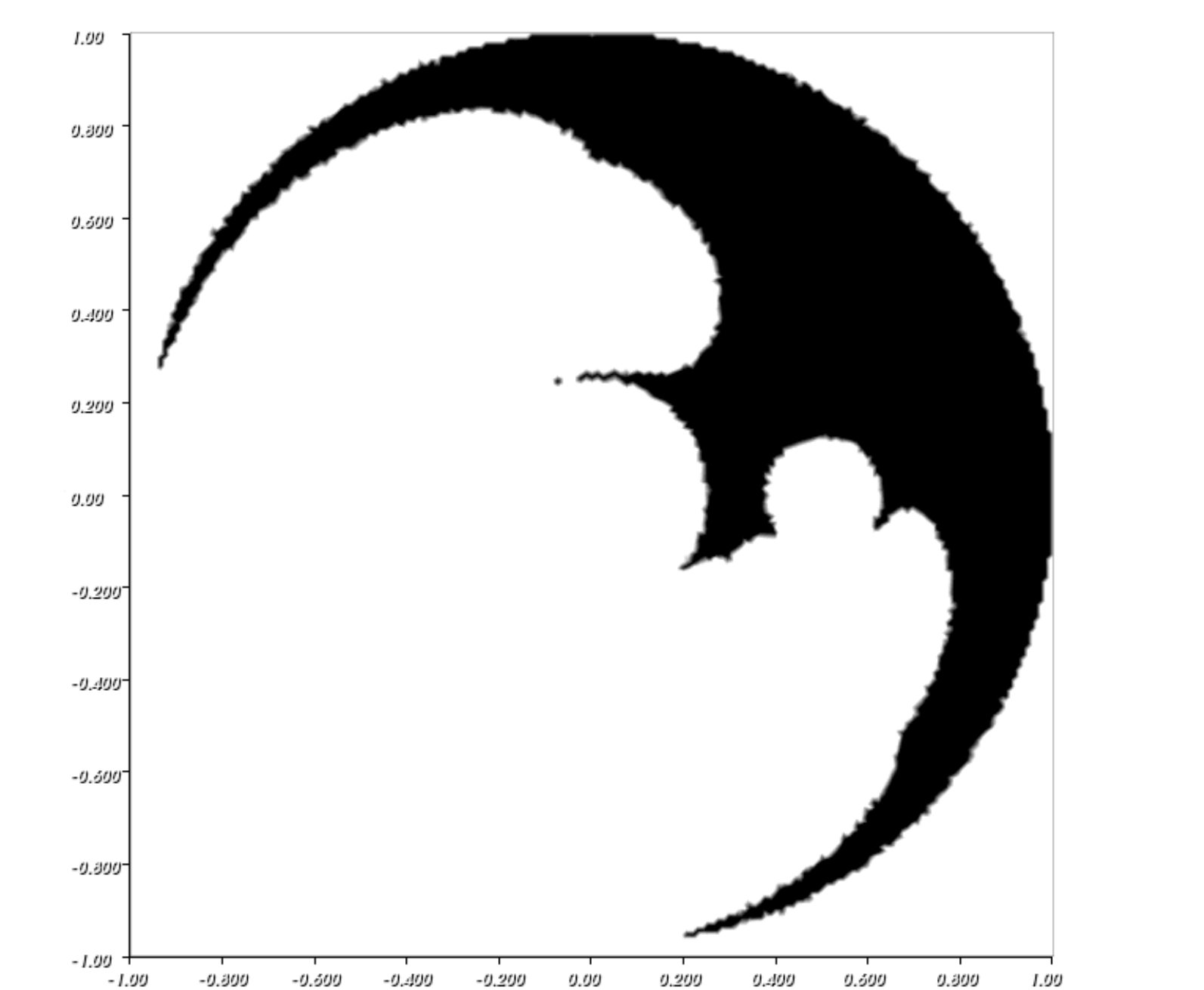}
\\
\includegraphics[width=5cm]{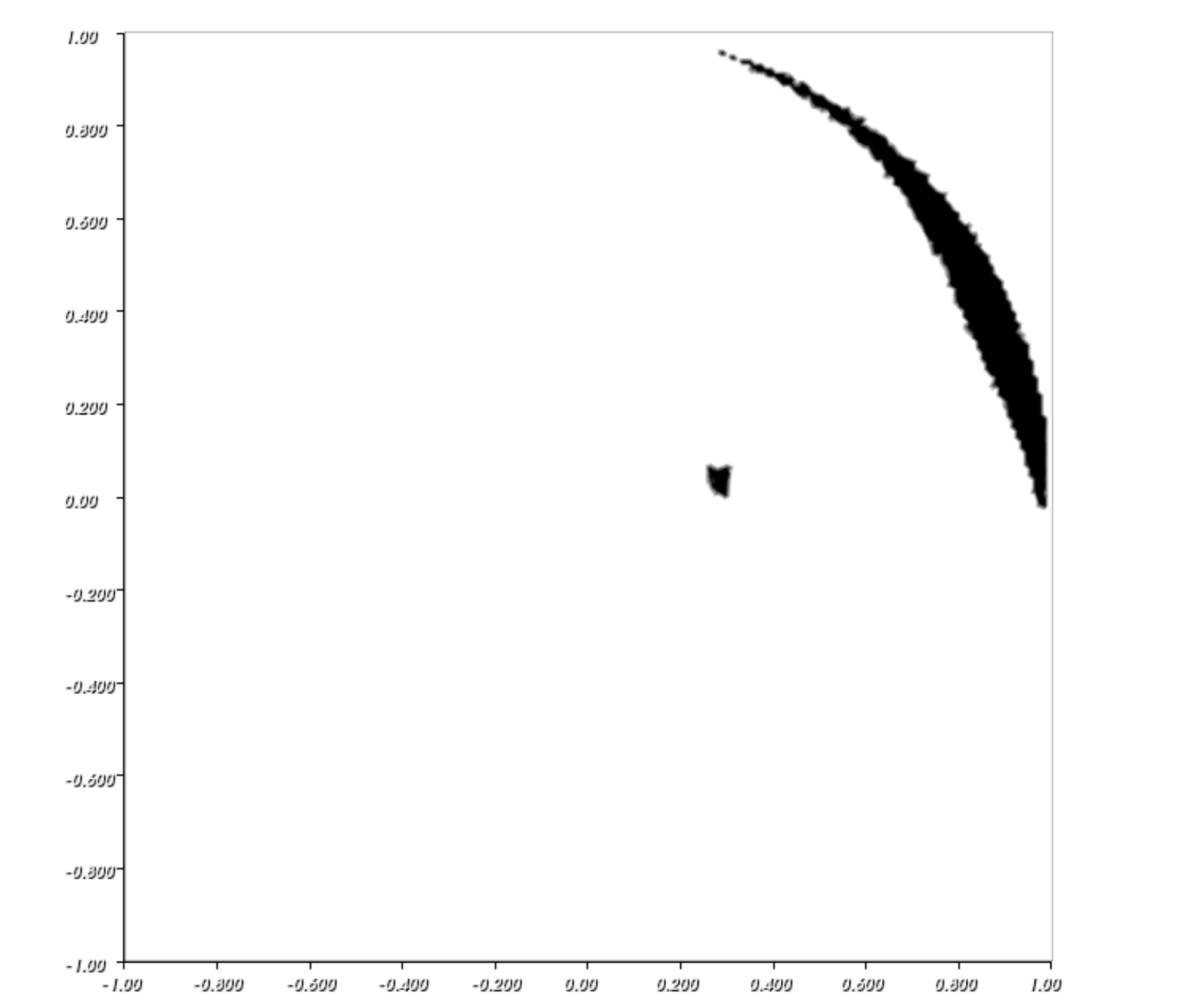}
\quad
\includegraphics[width=5cm]{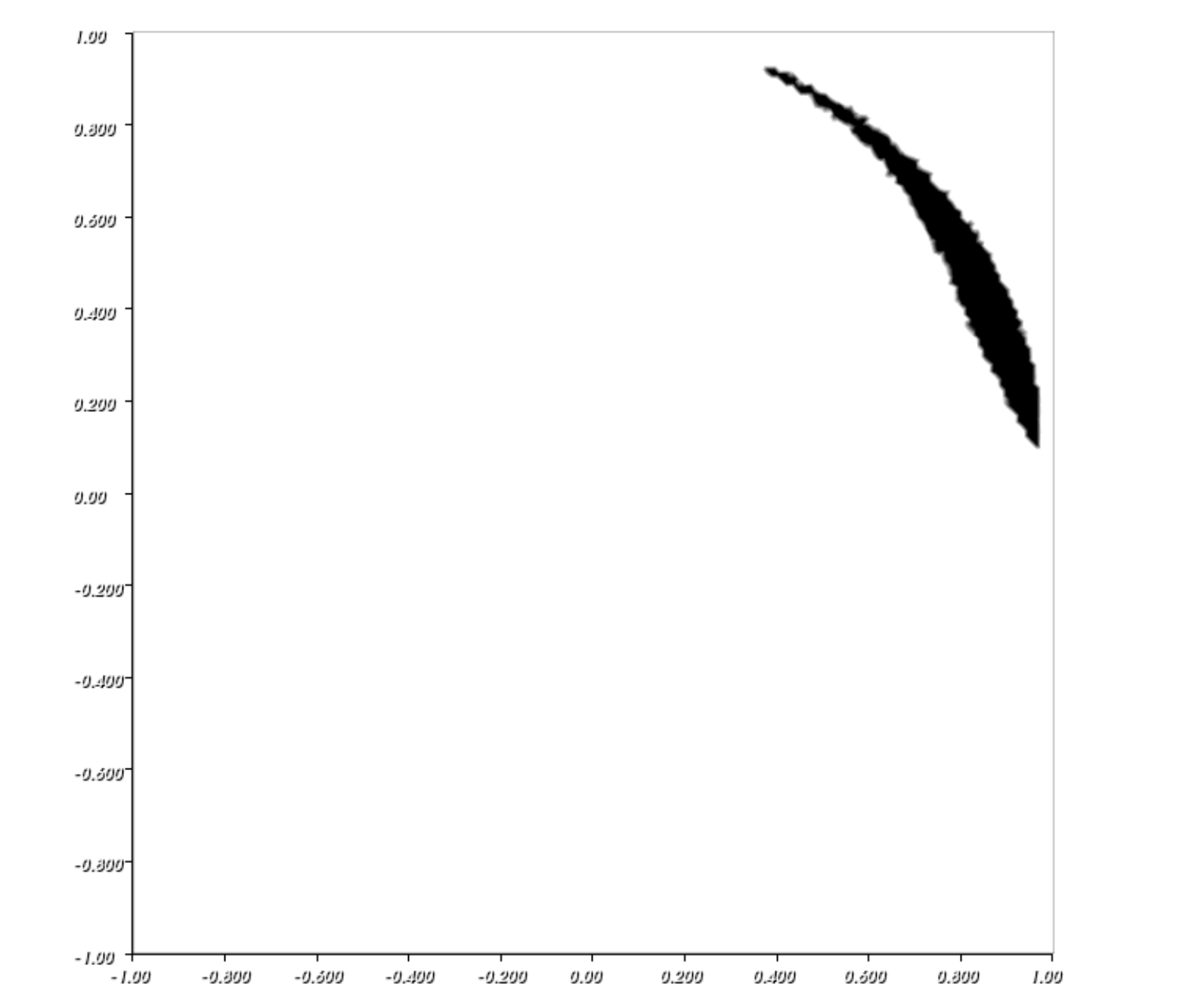}
\end{center}
\caption{Example 1. The initial domain with cost $2.29164$ (top, left) 
and intermediary domains in the line search
at the first iteration with cost $1.64609$ (top, middle),
respectively $0.133162$ (top, right); the domains $\Omega_n$ for $n=1,3$ (bottom) 
using the descent direction (\ref{4.2}).}
\label{fig:ex2_Omega}
\end{figure}

The cost function decreases rapidly at the first iterations
$j(g_0)=2.29164$, $j(g_1)=0.00083009$, $j(g_2)=0.000510025$, $j(g_3)=0.000379625$, but
for $n\geq 4$, $\Omega_n$ is similar to $\Omega_3$ and cost function decreases slowly
$j(g_8)=0.000171446$, $j(g_{11})=0.00012326$, $j(g_{14})=0.000100719$. The initial domain
and some computed domains  are presented in 
Figure \ref{fig:ex2_Omega}.

\newpage

As a second test, we use the descent direction
\begin{equation}\label{4.3}
d_n=R\left( \frac{1}{\epsilon} w_n\right)
\end{equation}
where $w_n$ is given by (\ref{4.2}) and $R:\mathbb{R}\rightarrow \mathbb{R}$ is
defined by
\begin{equation}\label{4.4}
R(r)=
\left\{
\begin{array}{ll}
-1+e^{r}, & r <0,\\
1-e^{-r}, & r \geq 0.
\end{array}
\right.
\end{equation}
The function $R$ is strictly increasing, $R\left(\mathbb{R}\right) =]-1,1[$, 
$R(-r)=-R(r)$ for all $r\in \mathbb{R}$. 

\begin{proposition}\label{prop:4.1} The direction $d_n$ defined by (\ref{4.3}) 
where $w_n$ is given by (\ref{4.2}) is a descent direction at $g_n$ for the cost function
$j(g)=\frac{1}{2}\int_\Omega (y_\epsilon-y_d)^2 d\mathbf{x}$.
\end{proposition}

\textbf{Proof.} 
The directional derivative of the cost function was introduced in the previous section.
In this particular case, the directional derivative of the cost function at $g_n$ 
in the direction $d_n$ is
$$
\int_D \left( H^\epsilon\right)^\prime(g_n)
\left[ 
\frac{1}{2}(y_\epsilon^n-y_d)^2
+\frac{1}{\epsilon}(y_\epsilon^n p^n + z_\epsilon^n q^n)
\right]d_n\, d\mathbf{x} .
$$
It can be rewritten as 
\begin{eqnarray*}
\int_D \left( H^\epsilon\right)^\prime(g_n) (-w_n)d_n\, d\mathbf{x}
& = & \int_D \left( H^\epsilon\right)^\prime(g_n) (-w_n) R\left( \frac{1}{\epsilon} w_n\right)
\, d\mathbf{x}\\
&=&-\epsilon \int_D \left( H^\epsilon\right)^\prime(g_n) 
\left( \frac{1}{\epsilon} w_n\right)
R\left( \frac{1}{\epsilon} w_n\right)\, d\mathbf{x} \leq 0.
\end{eqnarray*}
For the last inequality, we have used that $\left( H^\epsilon\right)^\prime(g_n) >0$ in $D$ and 
the property of the function $R$
$$
r R(r) \geq 0,\quad \forall r\in \mathbb{R}
$$
which gives $\left( \frac{1}{\epsilon} w_n\right)
R\left( \frac{1}{\epsilon} w_n\right)\geq 0$ in $D$.
Consequently, $d_n$ is a descent direction.
We remark that this derivative is zero, if and only if $w_n=0$ in $D$.\quad $\Box$

\medskip

\begin{figure}[ht]
\begin{center}
\includegraphics[width=5cm]{ex2_new_Omega0_axe.pdf}
\quad
\includegraphics[width=5cm]{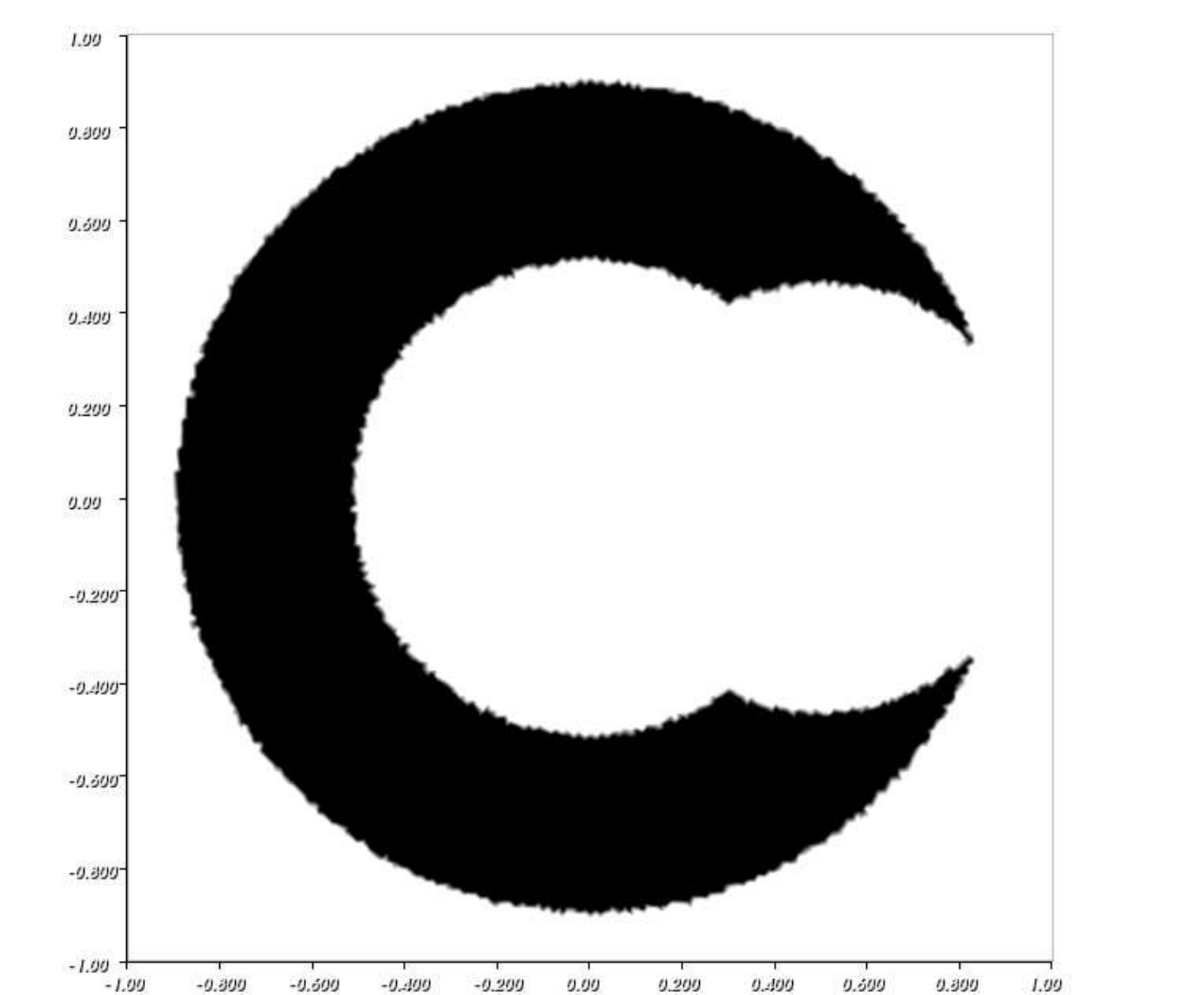}
\quad
\includegraphics[width=5cm]{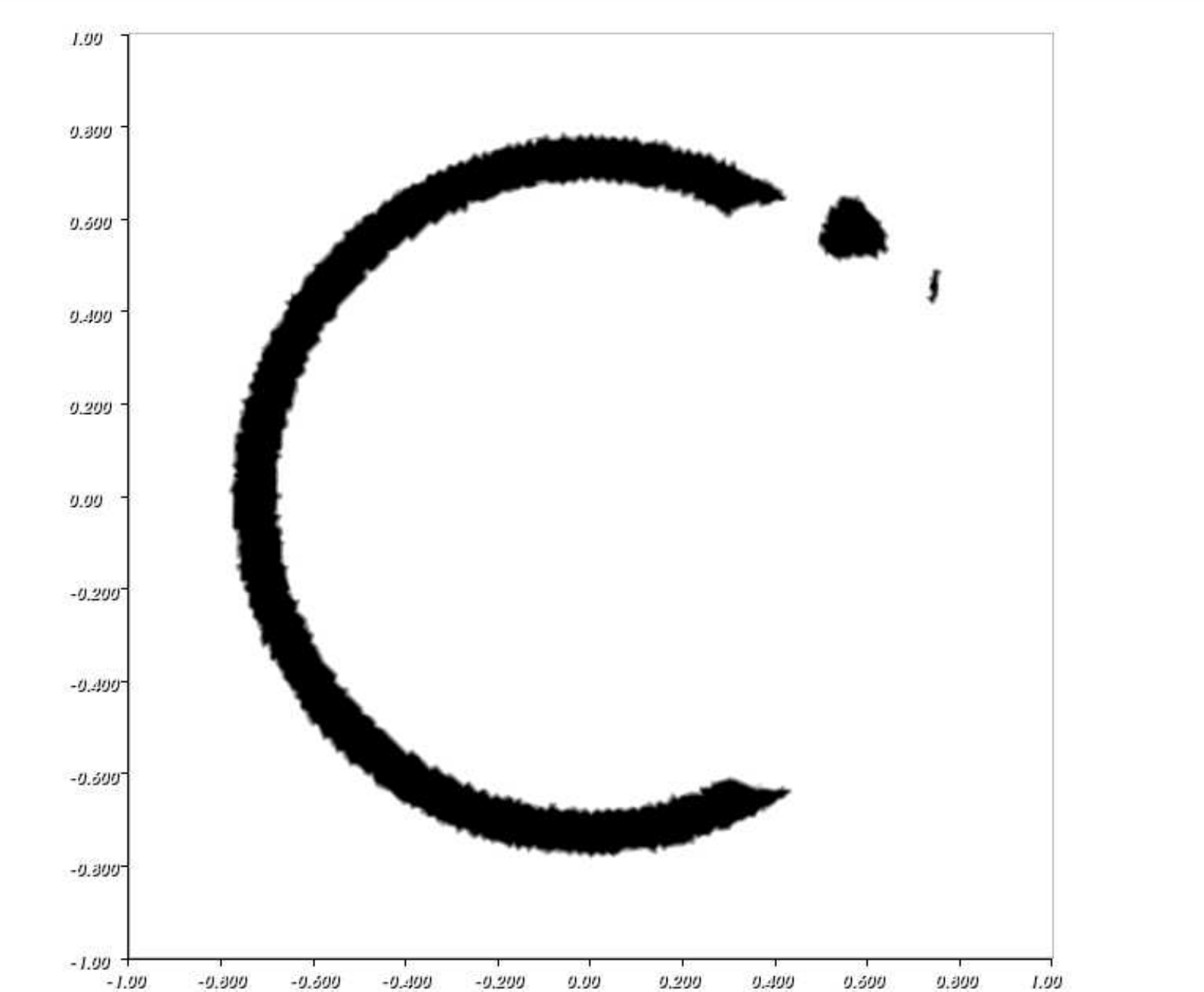}\\
\includegraphics[width=5cm]{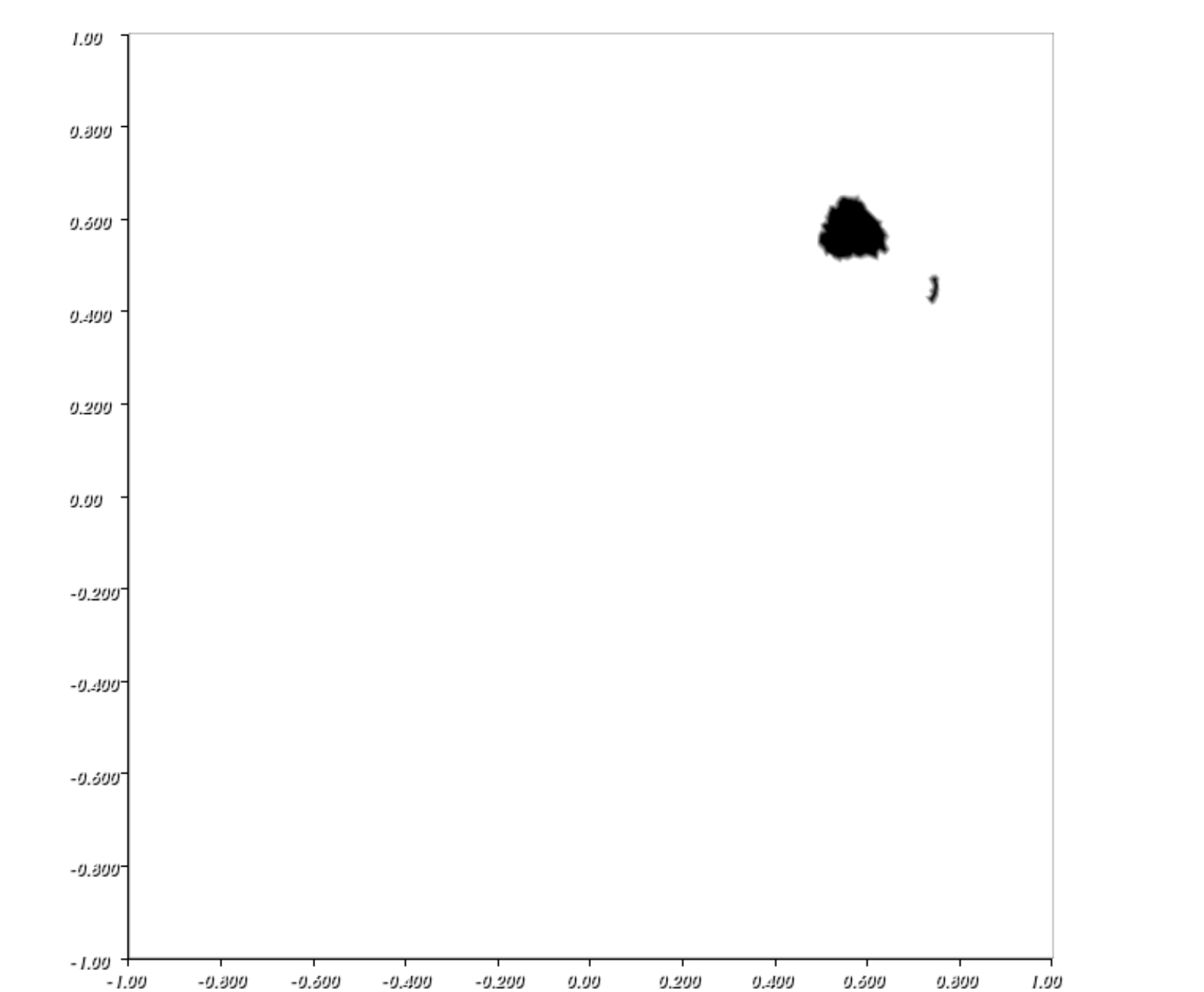}
\quad
\includegraphics[width=5cm]{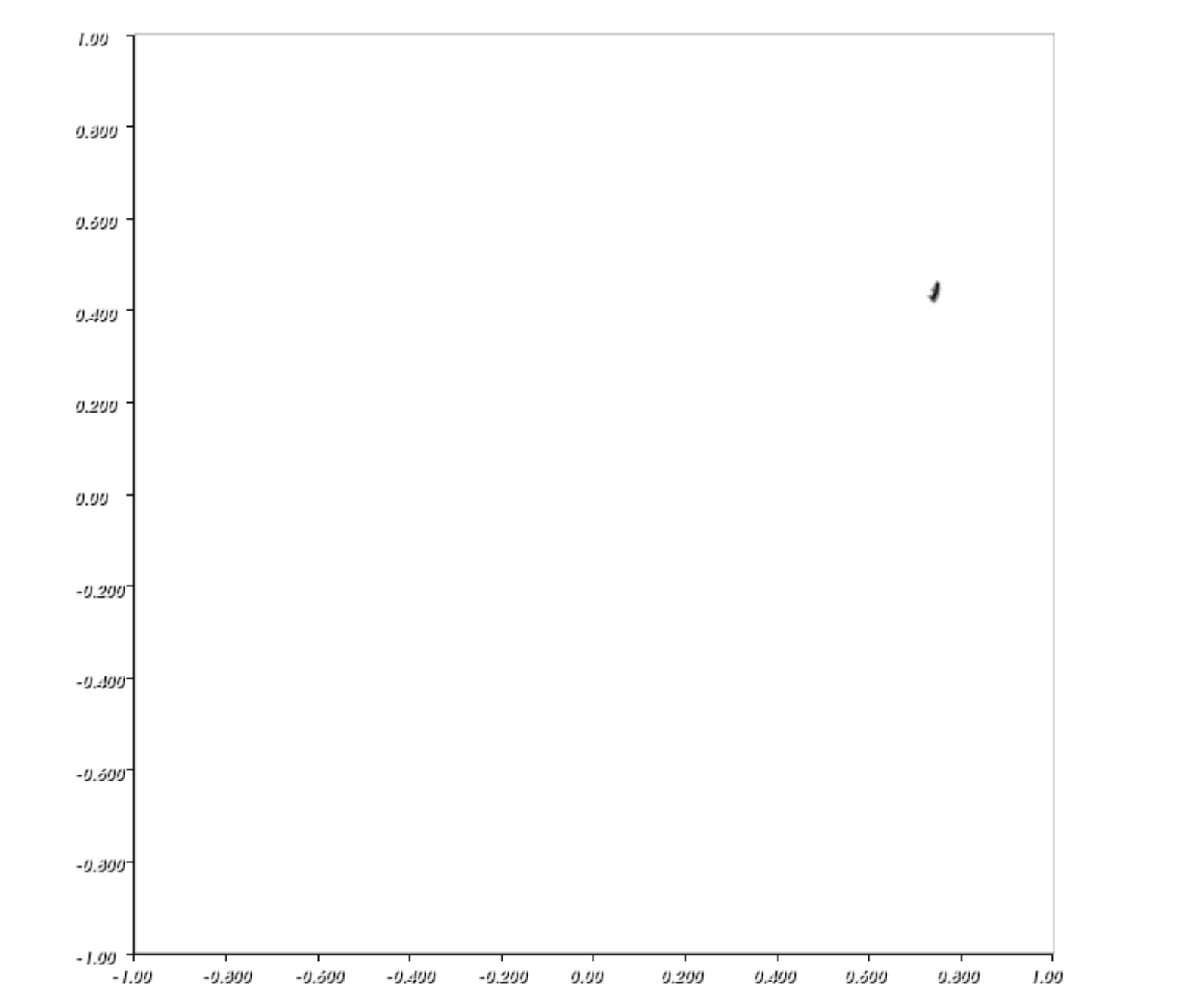}
\end{center}
\caption{Example 1. The domains $\Omega_n$ for $n=0,1,2$ (top) and $n=3,4$ (bottom)
using the descent direction (\ref{4.3}).}
\label{fig:ex2_R_Omega}
\end{figure}

In this second test, excepting the descent direction, the other parameters are the same as before.
The stopping test is obtained for $n=4$, the values of the cost function are:
$j(g_0)=2.29164$, $j(g_1)=1.23291$, $j(g_2)=0.295709$, $j(g_3)=1.66212e-05$,
$j(g_4)=5.0583e-11$. The computed domains are presented in Figure \ref{fig:ex2_R_Omega}. 
The optimal domain is the empty set and the optimal cost is zero, as obtained in both experiments.

\bigskip

\textbf{Example 2.}

We have again $D=]-1,1[\times ]-1,1[$. 
The load is $f=1$, the cost function is $j(g)=\int_\Omega (y_\epsilon-y_d) d\mathbf{x}$ where $y_d$ is given by
$$
y_d(x_1,x_2)=
\left\{
\begin{array}{rl}
 1, &  \hbox{if } \frac{1}{9} \leq x_1^2 + x_2^2 \leq \frac{1}{4} \\
-1, & \hbox{otherwise}. 
\end{array}
\right.
$$

We use for $D$ a mesh of 53360 triangles and 26981 vertices and 
for the approximation of $g$, $y$, $z$ we use piecewise linear 
finite element, globally continuous. 
The penalization parameter is $\epsilon=10^{-3}$. 

The cost functional is of type (\ref{3.8}) with 
$J\left(\mathbf{x},y_\epsilon(\mathbf{x})\right)=y_\epsilon-y_d$.
From Corollary \ref{cor:3.1} and Remark \ref{rem:3.1}, we get the
following descent direction 
\begin{equation}\label{4.5}
w_n=-\left[ 
(y_\epsilon^n-y_d)
+\frac{1}{\epsilon}(y_\epsilon^n p^n + z_\epsilon^n q^n)
\right].
\end{equation}

The sequence $\left( j(g_{n}) \right)_{n\in\mathbb{N}}$ is decreasing.
For the stopping test, we use: \texttt{if} $ j(g_{n+1}) > j(g_{n}) - tol$ 
\texttt{then} STOP, where $tol=10^{-6}$.

\begin{figure}[ht]
\begin{center}
\includegraphics[width=5cm]{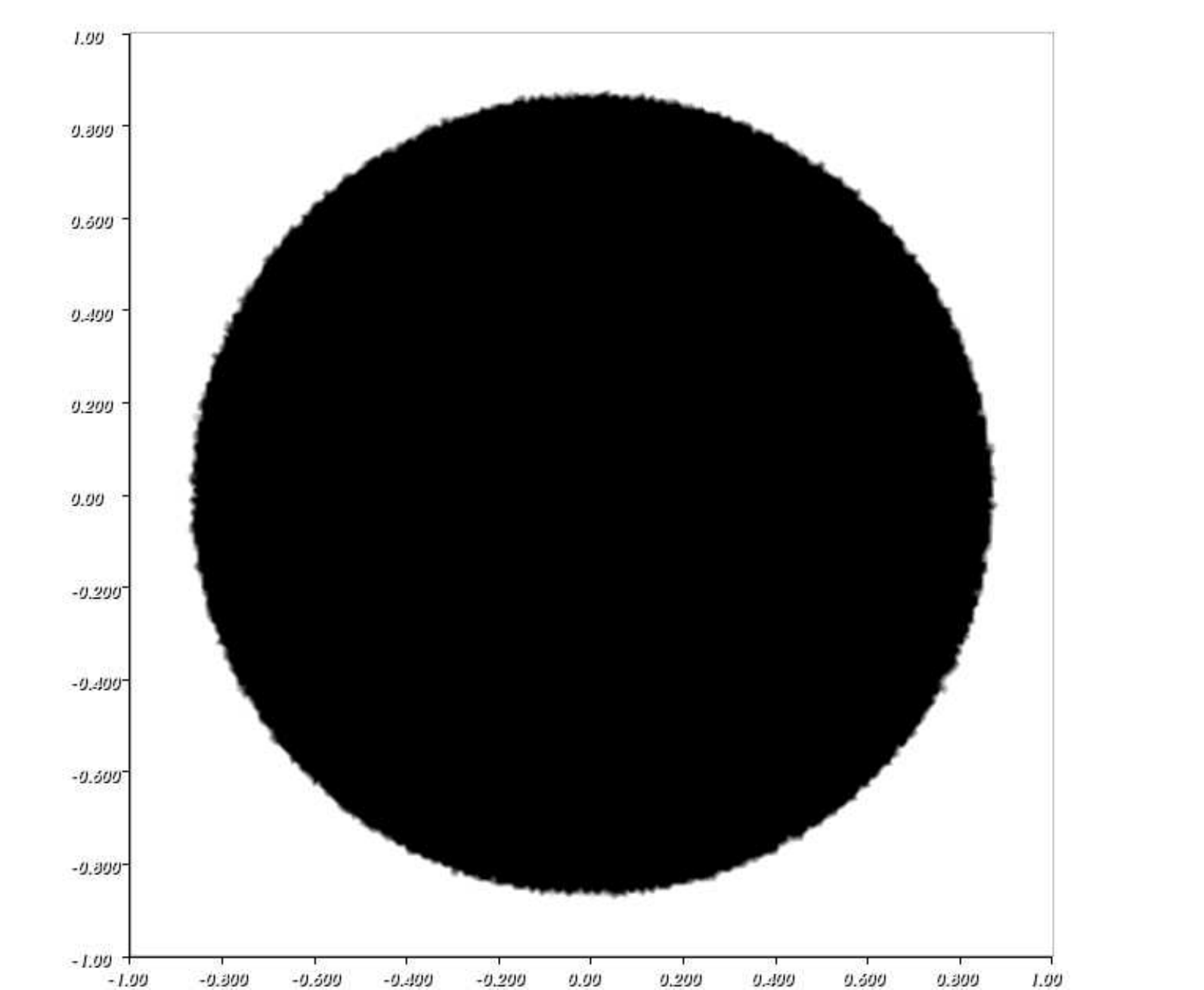}\quad
\includegraphics[width=5cm]{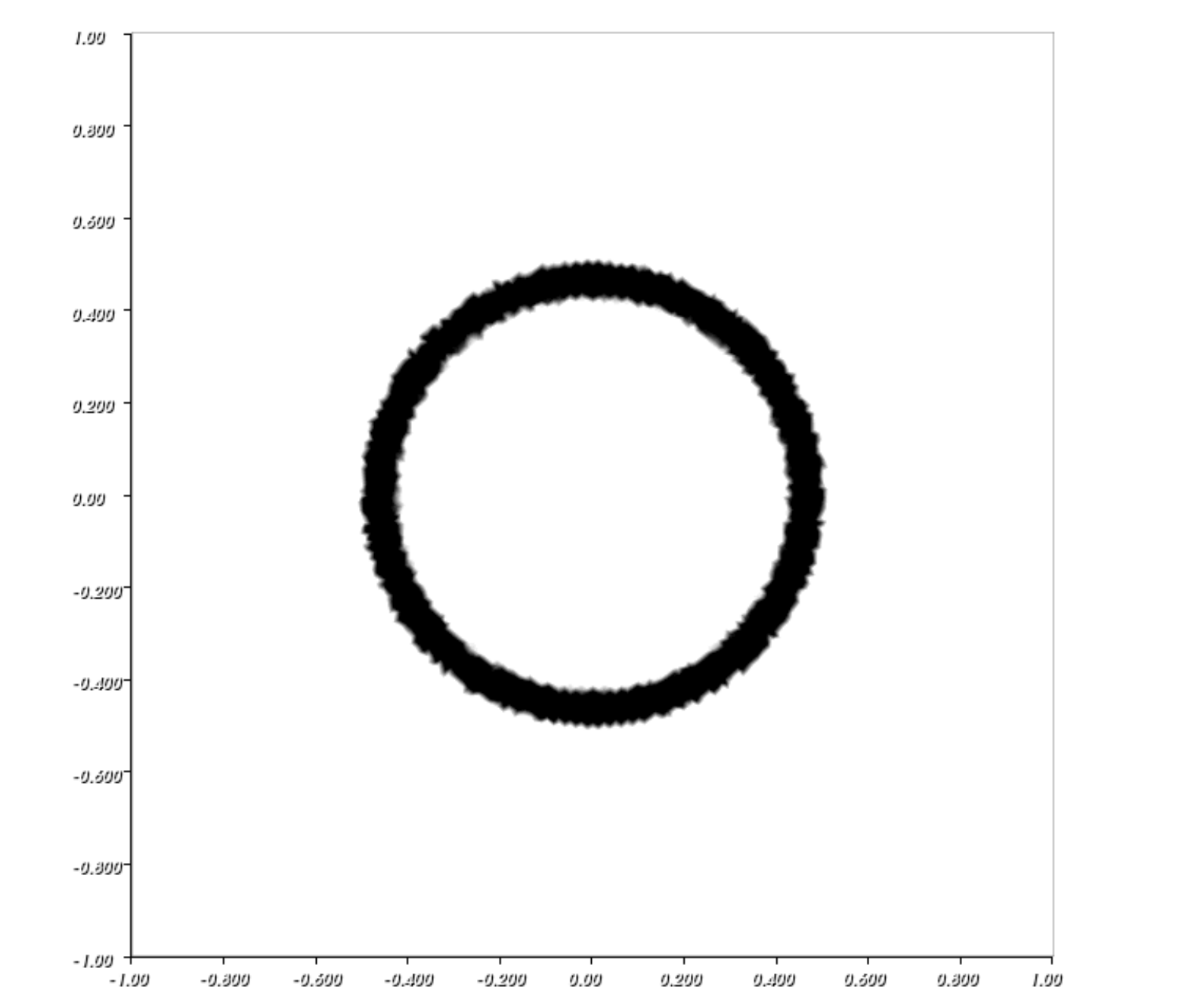}\quad
\includegraphics[width=5cm]{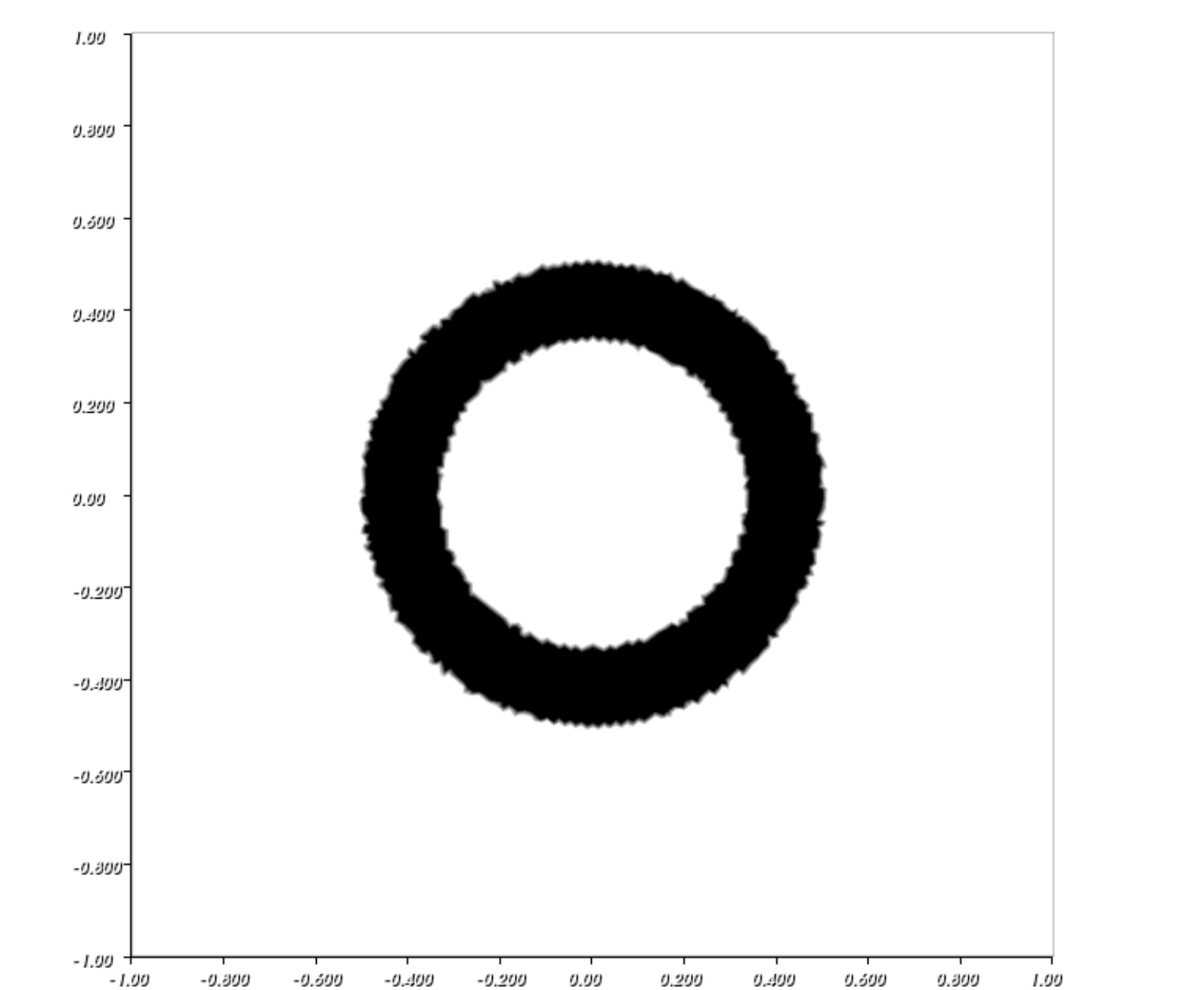}
\end{center}
\caption{Example 2. The initial domain $\Omega_0$ with cost 
$1.51761$ (left), intermediary domain in the line search with
cost  $-0.203754$ (middle) and optimal domain $\Omega_3$ with cost 
 $-0.423723$ (right), 
for initial parametrization $g_0(x_1,x_2)=-x_1^2-x_2^2+\frac{3}{4}$.}
\label{fig:iii_Omega}
\end{figure}

For the initial parametrization function $g_0(x_1,x_2)=-x_1^2-x_2^2+\frac{3}{4}$, that corresponds to a simply connected domain,
the stopping test is obtained for $n=3$, the values of the cost function are:
$j(g_0)=1.51761$, $j(g_1)=-0.417807$, $j(g_2)=-0.421269$, $j(g_3)=-0.423723$. 
Some computed domains are presented in Figure \ref{fig:iii_Omega}.

\newpage

\begin{figure}[ht]
\begin{center}
\includegraphics[width=5cm]{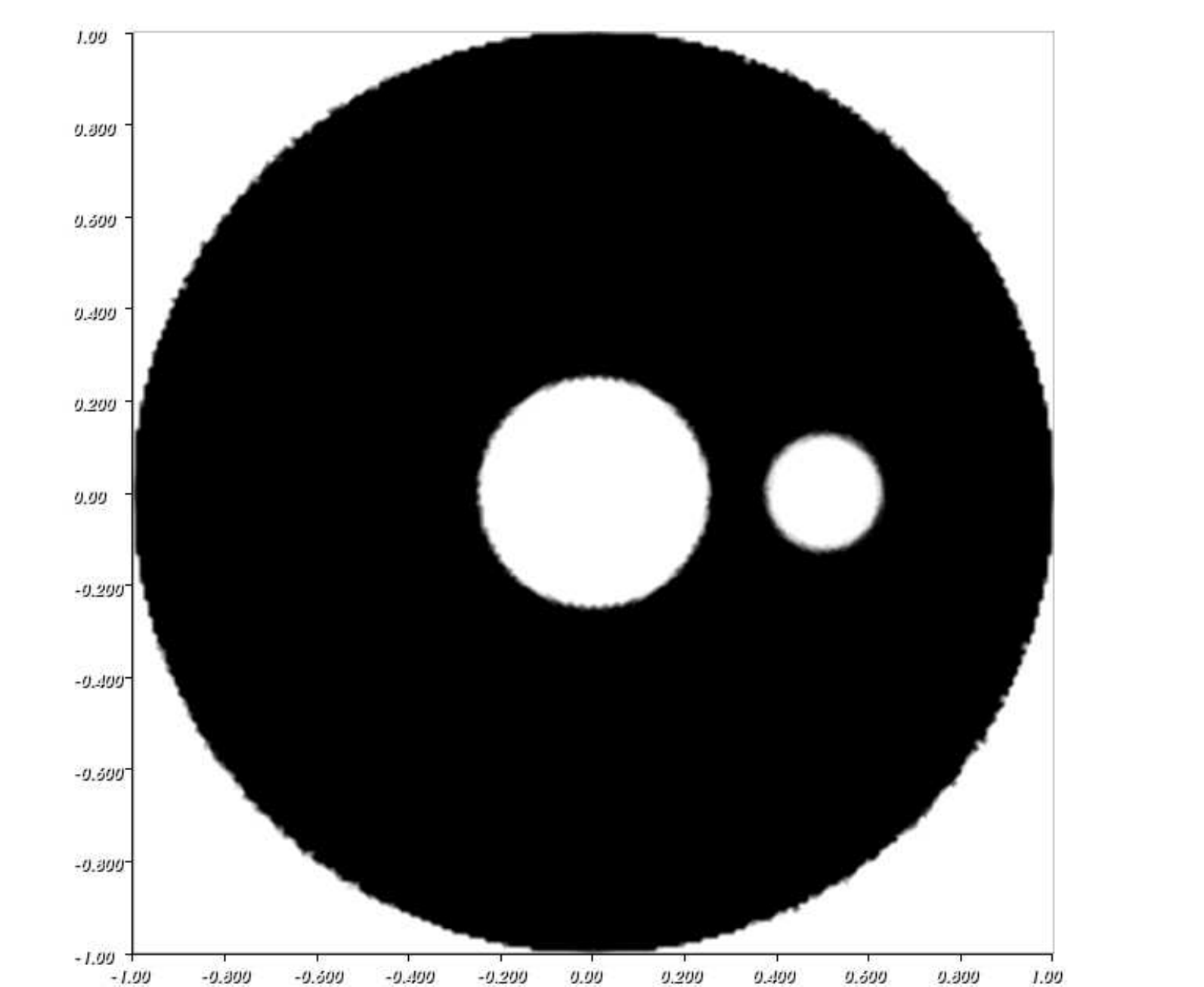}\quad
\includegraphics[width=5cm]{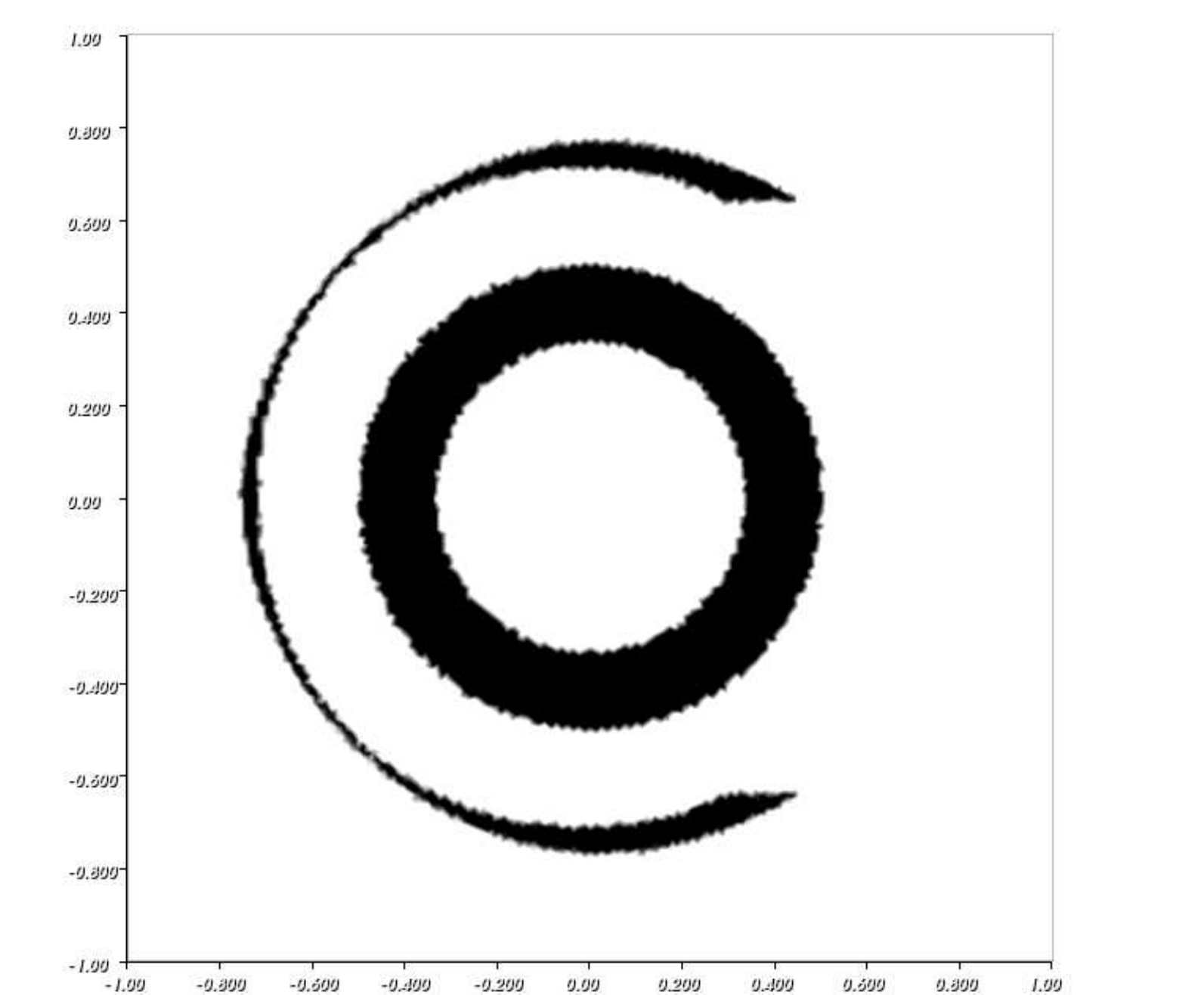}\quad
\includegraphics[width=5cm]{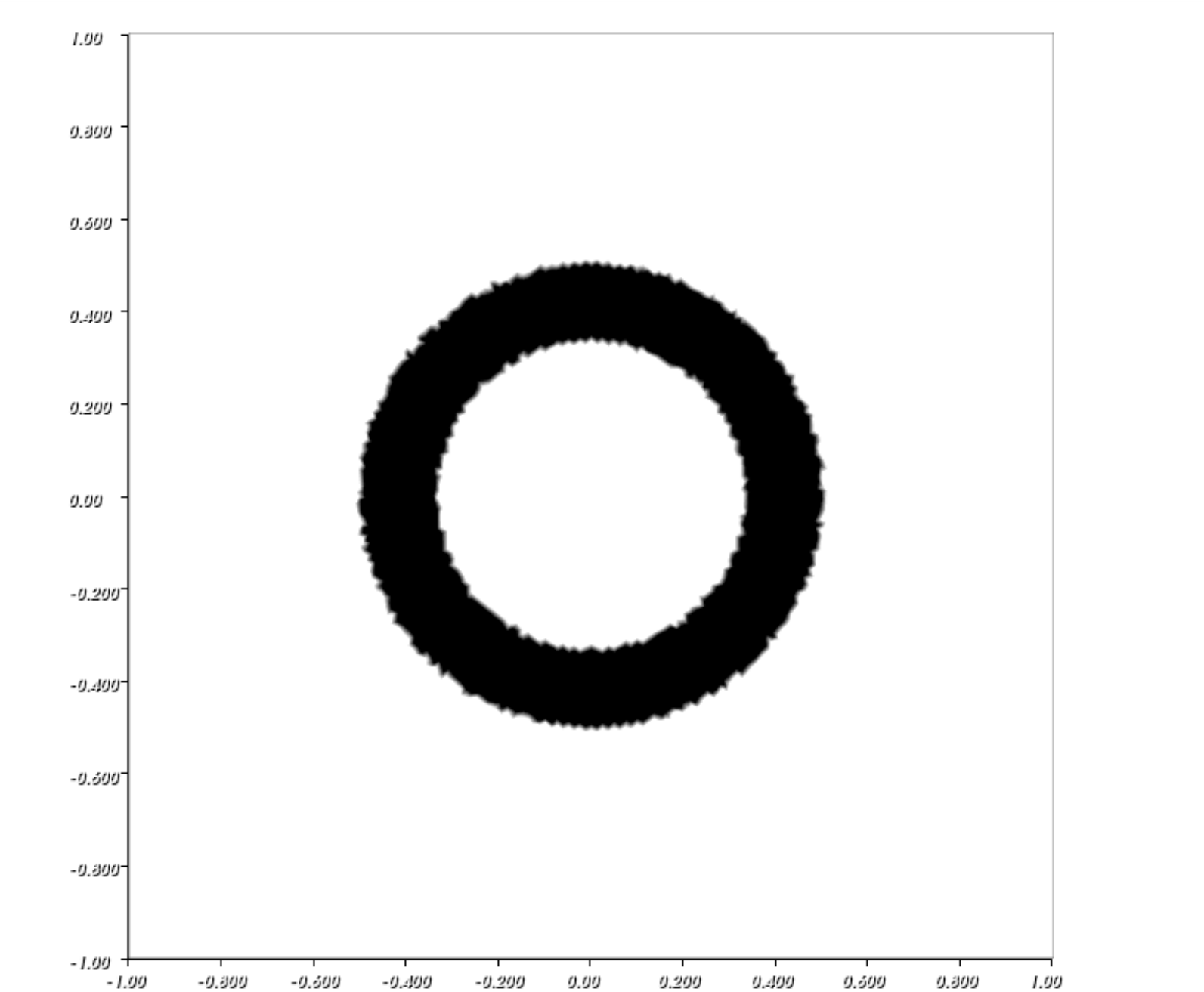}
\end{center}
\caption{Example 2. The initial domain $\Omega_0$ with cost 
$2.07908$ (left), intermediary domains with
cost  $-0.309447$ (middle) and $-0.424$ (right), 
for initial parametrization $g_0$ used in Example 1.}
\label{fig:iii_new_Omega}
\end{figure}

For the initial parametrization function  used in Example 1, 
the stopping test is obtained for $n=6$, the values of the cost function are:
$j(g_0)=2.07908$, $j(g_1)=-0.309447$, $j(g_2)=-0.424$, $j(g_3)=-0.424701$, 
$j(g_4)=-0.425225$, $j(g_5)=-0.425309$, $j(g_6)=-0.425331$. 
Some computed domains are presented in Figure \ref{fig:iii_new_Omega}.
The computed optimal cost depends slightly on $g_0$.

\bigskip

\textbf{Example 3.}

We have again $D=]-1,1[\times ]-1,1[$ and 
the cost function is $j(g)=\int_\Omega (y_\epsilon-y_d) d\mathbf{x}$.
 
The load is $f=2 \times 10^3$ and
$$
y_d(x_1,x_2)=
\left\{
\begin{array}{rl}
1, &  \hbox{if } x_1^2 + x_2^2 \leq \frac{1}{4}\\
-1, & \hbox{otherwise}. 
\end{array}
\right.
$$

We use the direction
\begin{equation}\label{4.6}
d_n=R\left( \frac{1}{\epsilon} w_n\right)
\end{equation}
where $w_n$ is given by (\ref{4.5}) and $R$ is defined by (\ref{4.4}).
As in Proposition \ref{prop:4.1}, it yields that $d_n$ is
a descent direction for the cost function $\int_\Omega (y_\epsilon-y_d) d\mathbf{x}$.
The other parameters are the same as in Example 1.

\begin{figure}[ht]
\begin{center}
\includegraphics[width=5cm]{iii_Omega0_axe.pdf}
\quad
\includegraphics[width=5cm]{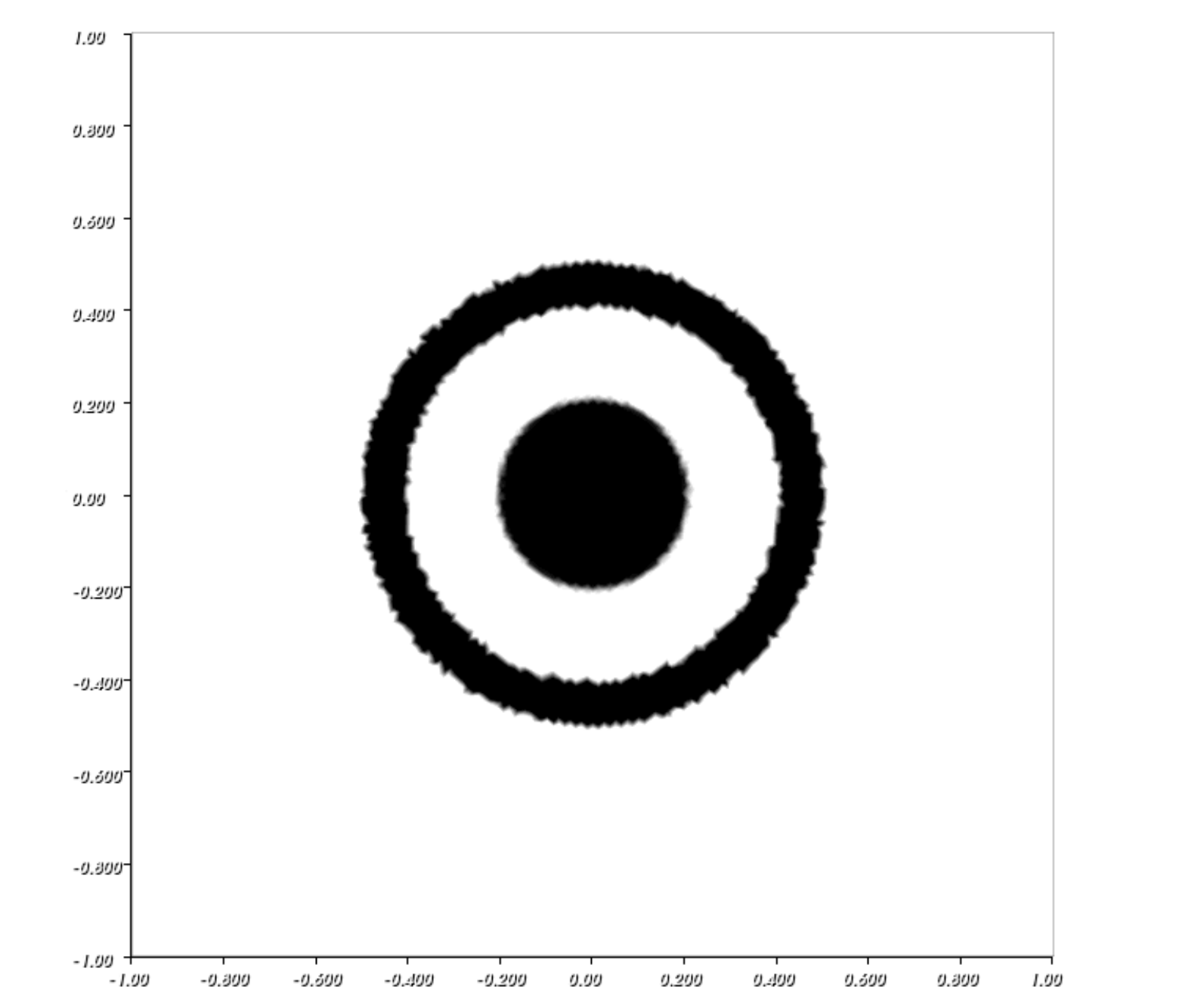}
\quad
\includegraphics[width=5cm]{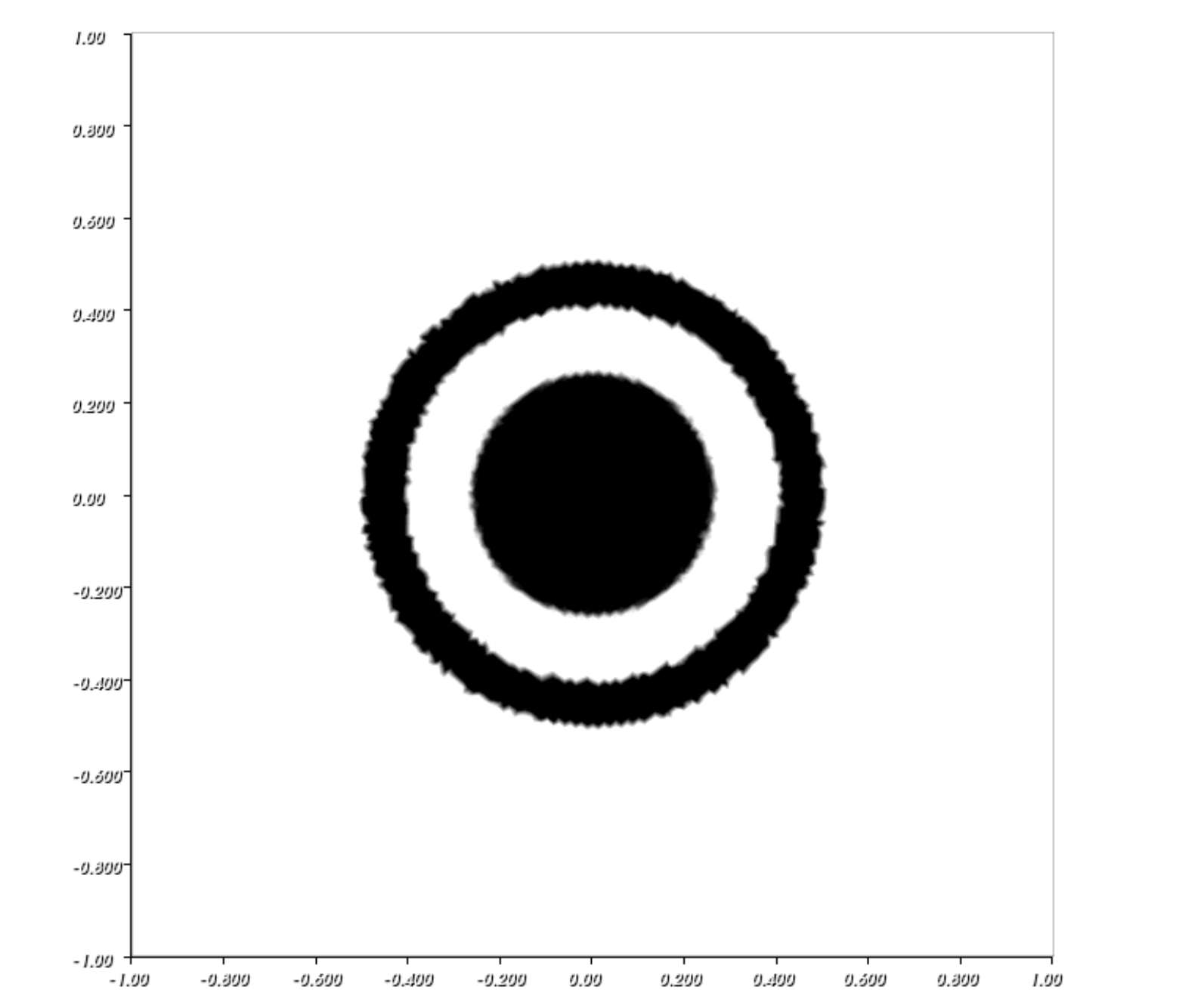}
\end{center}
\caption{Example 3. The domains $\Omega_n$ for $n=0,3,5$ using the descent direction (\ref{4.6})
and initial parametrization $g_0(x_1,x_2)=-x_1^2-x_2^2+\frac{3}{4}$.}
\label{fig:iv_Omega}
\end{figure}

For the initial parametrization function $g_0(x_1,x_2)=-x_1^2-x_2^2+\frac{3}{4}$, used in Example 2,
the stopping test is obtained for $n=5$, the values of the cost function are:
$j(g_0)=69.1791$, $j(g_1)=15.5425$, $j(g_2)=0.234407$, $j(g_3)=-0.34875$,
$j(g_4)=-0.385096$, $j(g_5)=-0.385548$. 
Some computed domains are presented in Figure \ref{fig:iv_Omega}.

\begin{figure}[ht]
\begin{center}
\includegraphics[width=5cm]{iii_new_Omega0_axe.pdf}
\quad
\includegraphics[width=5cm]{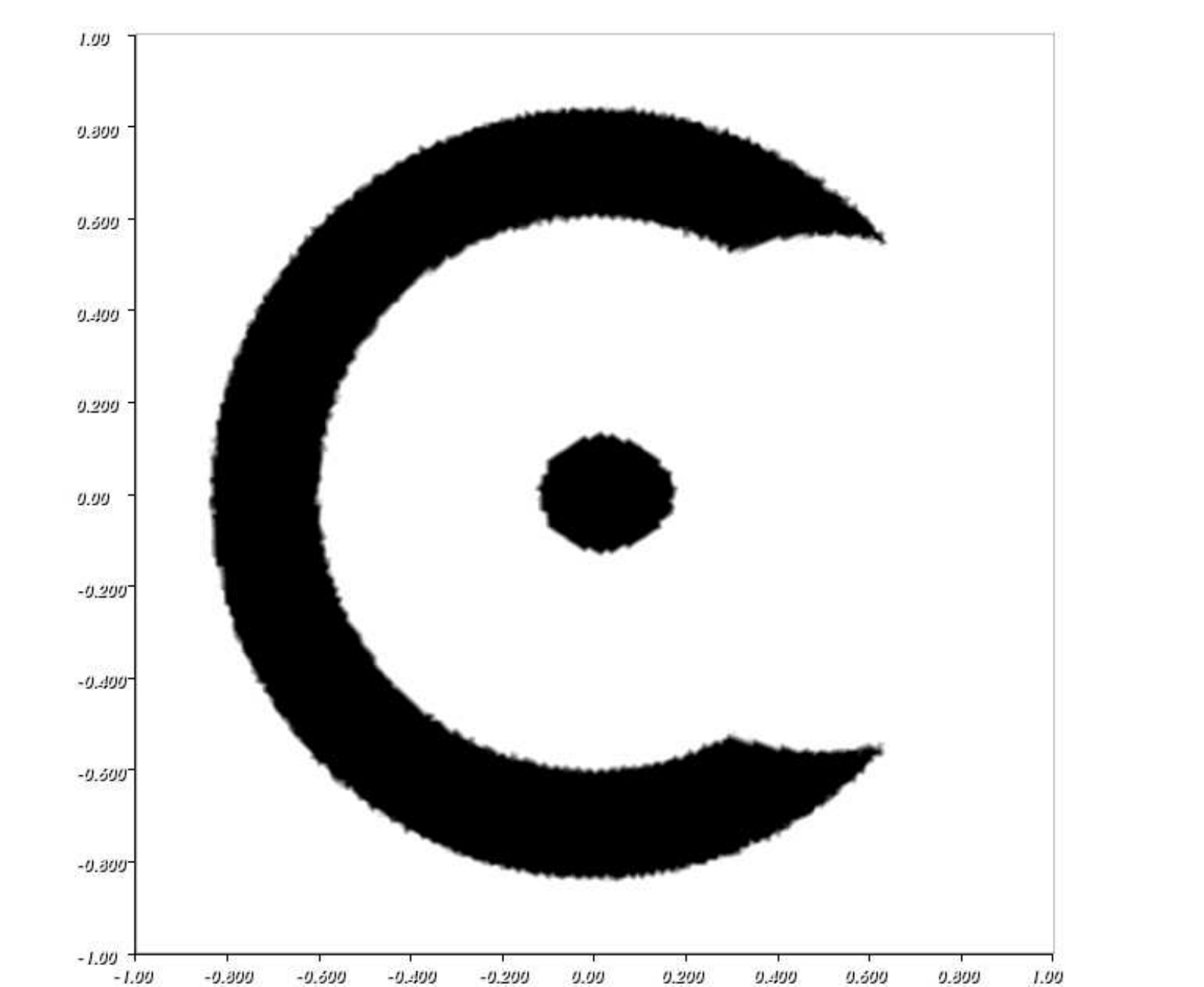}
\quad
\includegraphics[width=5cm]{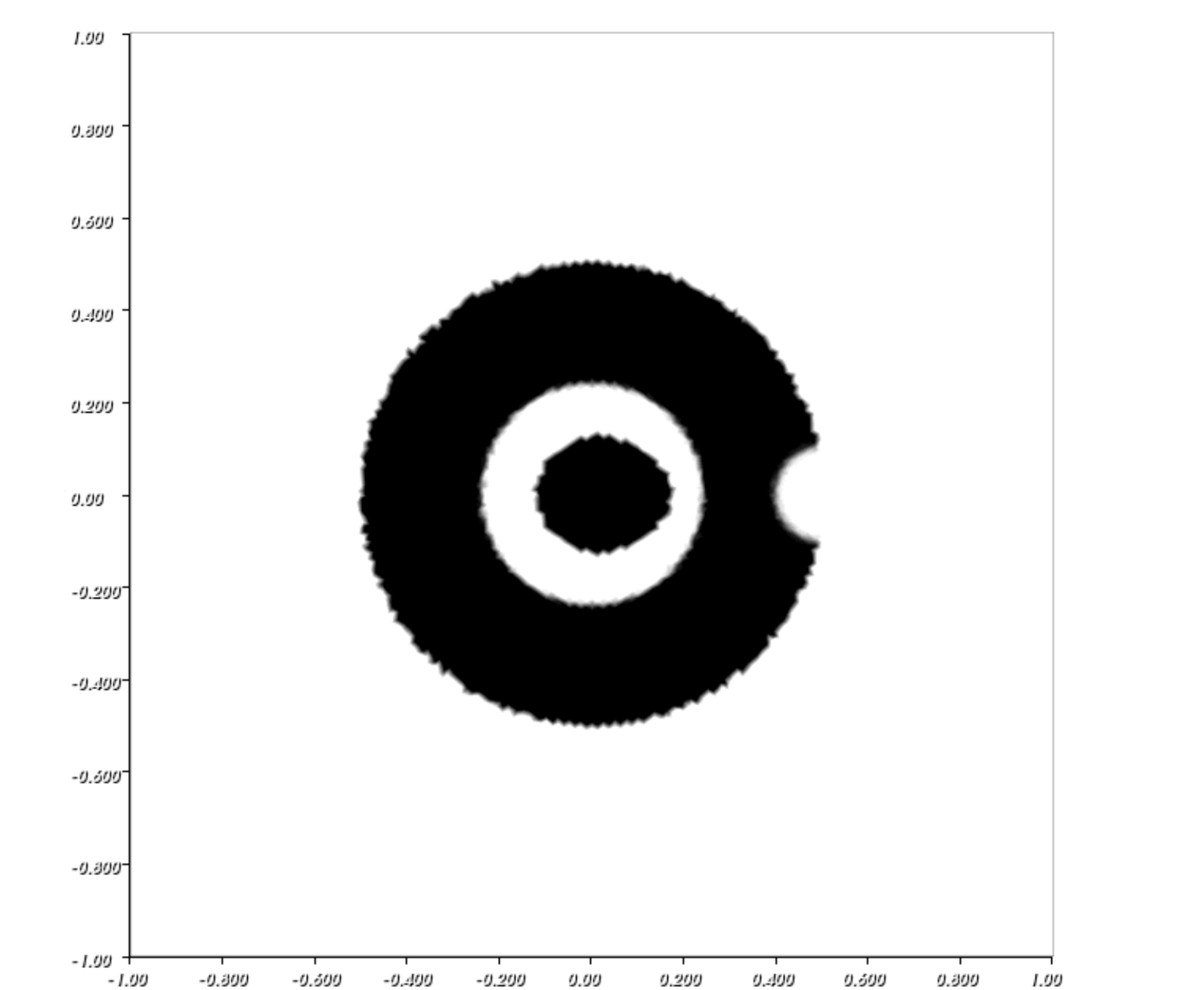}
\end{center}
\caption{Example 3. The domains $\Omega_n$ for $n=0,1,3$ using the descent direction (\ref{4.6})
for initial parametrization $g_0$ used in Example 1.}
\label{fig:iv_new_Omega}
\end{figure}

For the initial parametrization function used in Example 1, 
the stopping test is reached for $n=3$ and the values of the cost function are:
$j(g_0)=20.2385$, $j(g_1)=0.828123$, $j(g_2)=-0.509888$, $j(g_3)=-0.511685$. 
Some computed domains are presented in Figure \ref{fig:iv_new_Omega}.
We observe that the obtained result depends on $g_0$. Shape optimization problems are strongly non convex and ``local'' solutions are obtained, in general.

\section*{References}

\begin{thebibliography}{99}
\bibitem{Arnautu2000}
V. Arnautu, H. Langmach, J. Sprekels, D. Tiba, 
On the approximation and the optimization of plates. Numer. Funct. Anal. Optim. 
21 (2000), 
no. 3-4, 337--354. 

\bibitem{BST2012} 
M. Barboteu, M. Sofonea, D. Tiba, 
The control variational method for beams in contact with deformable obstacles, 
Z. angew. Math. Mech., 92 (2012) no. 1, pp. 25--40.

\bibitem{Chenais1975} 
D. Chenais, 
{On the existence of a solution in a domain identification 
problem}. {J. Math. Anal. Appl.} {52} (1975), no. 2, 189--219.

\bibitem{Delfour2001} 
M.C. Delfour, J.P. Zolesio, 
Shapes and Geometries, Analysis, Differential Calculus and Optimization, 
SIAM, Philadelphia, 2001.

\bibitem{Grisvard1985}
P. Grisvard, 
{Elliptic Problems in Nonsmooth Domains}. 
London, Pitman, 1985.

\bibitem{HMT2016} 
A. Halanay, C.M. Murea, D. Tiba, 
Existence of a steady flow of Stokes fluid past a linear elastic structure 
using fictitious domain, 
{J. Math. Fluid Mech.} 18 (2016) 397--413.

\bibitem{freefem++} F. Hecht,
New development in FreeFem++. 
{J. Numer. Math.} 20 (2012) 251--265.  
\texttt{http://www.freefem.org}

\bibitem{Kawohl1997}
B. Kawohl, J. Lang,
 Are some optimal shape problems convex? 
J. Convex Anal. 4 (1997), no. 2, 353--361.

\bibitem{Lions1969} 
J.-L. Lions, 
Quelques m\'ethodes de r\'esolution des probl\`emes aux limites non lin\'eaires,
Dunod, Paris, 
1969.

\bibitem{1992}
R. M\"akinen, P. Neittaanm\"aki, D. Tiba, 
On a fixed domain approach for a shape optimization problem. 
In: W. F. Ames, P. J. van der Houwen (Eds), Computational and Applied
Mathematics II: Differential Equations, North-Holland, Amsterdam, 1992, pp. 317--326.

\bibitem{Munoz2007}
J. Mu\~noz, P. Pedregal,
A review of an optimal design problem for a plate of variable thickness. 
SIAM J. Control Optim. 46 (2007), no. 1, 1--13.

\bibitem{MT2016}
C.M. Murea, D. Tiba,
A direct algorithm in some free boundary problems. 
{J. Numer. Math.} 24 (2016) 253--271.

\bibitem{Tiba2009} 
P. Neittaanm\"aki, A. Pennanen, D. Tiba, 
{Fixed domain approaches in shape optimization problems with Dirichlet 
boundary conditions},  
{Inverse Problems} 25 (2009) 1--18.

\bibitem{Tiba2006} 
P. Neittaanm\"aki, J. Sprekels, D. Tiba, 
{Optimization of elliptic systems. Theory and applications}, Springer, 
New York, 2006.

\bibitem{Tiba2012} 
P. Neittaanm\"aki, D. Tiba, 
Fixed domain approaches in shape optimization problems, 
{Inverse Problems} 28 (2012) 1--35.

\bibitem{Pironneau1984} 
O. Pironneau, 
{Optimal shape design for elliptic systems}, Springer, Berlin, 1984.

\bibitem{ptiba}
P. Philip, D. Tiba,
Shape optimization via control of a shape function on a fixed domain: 
theory and numerical results. In: S. Repin, T. Tiihonen, T. Tuovinen (Eds), 
Numerical methods for differential equations, optimization and technological problems, 
Computational methods in applied sciences 27, Springer Verlag, Dordrecht, 2013,
pp. 305--320.


\bibitem{ST} 
M. Sofonea, D. Tiba,
The control variational method for contact of Euler-Bernoulli beams, 
{Bull. Transilvania Univ. Bra\c{s}ov} vol. 2 (51), Series
III (2009), p.127--136.


\bibitem{Sokolowski1992}
J. Sokolowski, J.P. Zolesio, 
{Introduction to Shape Optimization. Shape Sensitivity Analysis}, 
Springer, Berlin,
1992.

\bibitem{Sprekels2003}
J. Sprekels, D. Tiba, 
Optimization of clamped plates with discontinuous thickness. Optimization and control 
of distributed systems. {Systems Control Lett.} 48 (2003), no. 3-4, 289--295.

\bibitem{tiba2013}
D. Tiba,
Domains of class C: properties and application. 
{Ann. of the Univ. of Bucharest (Ser. Math.)}, 4(LXII) (2013), 89--102.

\bibitem{2017}
D. Tiba,
Neumann boundary conditions in shape optimization,
Pure Appl. Funct. Anal.,
accepted (2017)

\end{thebibliography}
\end{document}